\documentclass[a4paper,11pt]{amsbook}

\usepackage[latin1]{inputenc}
\usepackage[english]{babel}
\usepackage[pdftex]{hyperref}
\usepackage{a4wide}
\usepackage{graphicx}
\usepackage{amsmath,amssymb,amsthm, amsfonts}
\usepackage{fancyhdr}
\hypersetup{colorlinks 	= true,
			linkcolor 	= black,
			anchorcolor = black,
    		citecolor 	= blue,
    		filecolor 	= black,
    		pagecolor 	= black,
    		urlcolor 	= black}

\renewcommand{\Re}{\mathrm{Re}\ }
\renewcommand{\Im}{\mathrm{Im}\ }
\newtheorem{thm}{Theorem}[chapter]
\newtheorem{lemm}{Lemma}[chapter]
\newtheorem{cor}{Corollary}[chapter]
\theoremstyle{definition}

\theoremstyle{remark}
\newtheorem{rem}{Remark}

\newcommand{\LMUTitle}[9]{
  \thispagestyle{empty}
  \vspace*{\stretch{1}}
  {\parindent0cm
   \rule{\linewidth}{.7ex}}
  \begin{flushright}
  	\vspace*{\stretch{1}}
   	\sffamily\bfseries\LARGE #1\\
   	\vspace*{\stretch{1}}
   	\sffamily\bfseries\large
   	#2
   	\vspace*{\stretch{1}}
  \end{flushright} 
  \rule{\linewidth}{.7ex}
  \vspace*{\stretch{5}}
  \begin{center}
   \includegraphics[width=2in]{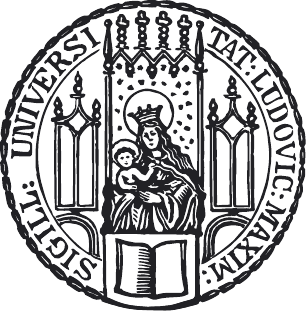}
  \end{center}
  \vspace*{\stretch{1}}
  \begin{center}\sffamily\huge{Diplomarbeit}\end{center}
  \begin{center}\sffamily\LARGE{#5}\end{center}
  \cleardoublepage
}

\renewcommand{\chaptermark}[1]%
         {\markleft{#1}}

\renewcommand{\sectionmark}[1]%
         {\markright{\thechapter.\thesection\ #1}}

\numberwithin{equation}{chapter}
\numberwithin{figure}{chapter}

\begin{document}

\frontmatter
\LMUTitle
      {On the Existence Theory of Hilbert Space Valued Diffusion Processes}		
      {Günter Hinrichs}                   
      {München}                   	
      {Mathematischen Institut}		    
      {M\"unchen 2010}                 	
      {26. Juli 2010}                   
      {Prof. Dr. Detlef D\"urr}            	
      {Prof. Dr. Vitali Wachtel}              	

\thispagestyle{empty}
\begin{center}
\LARGE{\textbf{On the Existence Theory of Hilbert Space Valued Diffusion Processes}}
\end{center}
\begin{center}
\Large{Günter Hinrichs}
\end{center}
\begin{center}
\today
\end{center}
\vspace{2cm}
\begin{center}
\begin{minipage}{0.8\textwidth}
\begin{center}
\textbf{Zusammenfassung}
\end{center}
 Es wird gezeigt, dass man die Lösung gewisser linearer stochastischer Differentialgleichungen in Hilberträumen, und zwar solcher mit beschränkten Operatoren sowie der in \cite{Holevo} und \cite{Mora} betrachteten konservativen stochastischen Schrödingergleichung, analog zur Lie-Trotter-Produktformel aus der Halbgruppentheorie gewinnen kann, indem man die Gleichungen in einen ``deterministischen'' und einen ``stochastischen'' Anteil aufspaltet und die zugehörigen Lösungsflüsse abwechselnd iterativ auf die gewünschte Anfangsbedingung anwendet.

\textbf{Betreuer:} Prof. Dr. Detlef Dürr
\end{minipage}
\end{center}

\vspace{2cm}

\begin{center}
\begin{minipage}{0.8\textwidth}
\begin{center}
\textbf{Abstract}
\end{center}
 We prove that the solution of certain linear stochastic differential equations in Hilbert spaces, namely those with bounded operators as well as the conservative stochastic Schrödinger equations analyzed in \cite{Holevo} and \cite{Mora}, can be obtained - along the lines of the Lie-Trotter product formula from semigroup theory - by splitting the equation into a ``deterministic'' and a ``stochastic'' part and alternately applying the corresponding solution flows in an iterative manner to the initial value.

\textbf{Advisor:} Prof. Dr. Detlef Dürr
\end{minipage}
\end{center}

\cleardoublepage

\tableofcontents
\cleardoublepage

\mainmatter
\setcounter{page}{1}

\chapter{Introduction}

The main subject of this work is the stochastic Schrödinger equation
\begin{equation}\label{346}
 \textnormal{d}\psi_t = (-iH-\frac{1}{2}A²)\psi_t\textnormal{d}t + A\psi_t\textnormal{d}\xi_t .
\end{equation}
for a Hilbert space valued diffusion process $\psi_t$ with standard Wiener process $\xi_t$ and, in general, unbounded self-adjoint operators H and A. This is an extension of the usual Schrödinger equation
\begin{equation}\label{347}
\textnormal d\psi_t = -iH\psi_t \textnormal dt
\end{equation}
by a stochastic term that will make the wavefunction ``collapse''.
Note that solving the isolated ``stochastic part''
\begin{equation}\label{348}
\textnormal d\psi_t=A\psi_t\textnormal d\xi_t - \frac{1}{2}A²\psi_t\textnormal dt
\end{equation}
of this equation is considerably simpler than solving the whole equation - in an abstract framework as well as if one wants to find an explicit solution formula. The same is true for the Hamiltonian part \eqref{347}.
Therefore I wondered whether one can find an alternative access to the solution of the whole equation by constructing it out of these two partial solutions. Such a method would reflect the intuitive notion leading to equations like \eqref{346} that the two mechanisms \eqref{347} and \eqref{348}, when combined additively, act simultaneously on $\psi_t$.

Let me remind you that the analogous question for deterministic equations of the form
\begin{equation}\label{784}
 \dot x(t)=(A+B)x(t)
\end{equation}
is settled by the Lie-Trotter product formula, which, in the presence of unbounded operators, holds true e.g. in the following form (see \cite{Chernoff}, \cite{Reed1}):
\begin{thm}
 Let A, B and $C:=\overline{A+B}$ be generators of contraction semigroups, noted as $e^{tA}$ etc., on a Banach space.Then $s-\lim_{n\to\infty}(e^{\frac{t}{n}A}e^{\frac{t}{n}B})^n=e^{tC}$ holds true with ``s-lim'' denoting the strong operator limit.
\end{thm}
Roughly speaking, this means that one gets the solution of \eqref{784} by subjecting the initial value alternately for infinitesimal time intervals to the solution flows $e^{tA}$ and $e^{tB}$ which belong to the single operators.
Chernoff has proven the Trotter formula in a more general framework, replacing $F(t):=e^{\frac{t}{n}A}e^{\frac{t}{n}B}$ by an arbitrary contraction valued mapping that fulfills $F(t)=id+tC+o(t)$ for $t\to0$:
\begin{thm}
 Let F be a strongly continuous mapping from $[0,\infty)$ to the set of contraction operators on a Banach space such that $\overline{F'(0)}$, the clousre of its strong derivative at 0, generates a contraction semigroup. Then $s-\lim_{n\to\infty}F\left(\frac{t}{n}\right)^n=e^{t\overline{F'(0)}}$ .
\end{thm}
Besides the classical Trotter formula, this theorem comprises e.g. the cases $F(t)=(id-tC)^{-1}$ and $F(t)=[(id-tA)(id-tB)]^{-1}$ . 

Our goal is to prove similar results for linear stochastic differential equations (SDE). In chapter \ref{ch:Loesung}, after having taken a look at the solution of the isolated stochastic part of \eqref{346}, we develop the theory of weak solutions of linear dissipative SDE - a class of equations comprising \eqref{346} - going back to \cite{Holevo} and mention the results of Mora und Rebolledo (\cite{Mora}) about their regularity. Chapter \ref{ch:Produkt} begins with a generalization of Chernoff's theorem to linear SDE which have matrices or bounded operators as coefficients. A result of this kind does not seem to have been published so far. Then we combine our proof with ideas form \cite{Holevo} and \cite{Mora} concerning the handling of unbounded operators and arrive at a product formula for the stochastic Schrödinger equation \eqref{346}.

Shortly before we finished this work, Smolyanov et al. have republished the article \cite{Gough}, where they also find a product formula for stochastic Schrödinger equations. However, they admit only Hamiltonians arising as the Weyl quantization of a classical phase space function and bounded collapse operaors instead of our unbounded operator A.

In chapter \ref{ch:Beispiel} we apply our product formula to two prototypes of quantum mechanical collapse models, namely the time-continuous QMUPL model (\cite{Diosi2}) and the discrete GRW model (\cite{Ghirardi1}), and shows that the first can be constructed as a suitable continuum limit of the latter.

\chapter{Solution theory for linear SDE with unbounded operators}\label{ch:Loesung}

\section{SDE with one operator}\label{sec:Einop}

Let $\xi_t$ be a Brownian motion on a probability space $(\Omega, \mathfrak{F}, \mathbb{P})$. A simple example for an infinite dimensional SDE is the equation
\begin{equation}\label{geom}
\begin{split}
 \mathrm{d}f_t(x) = & xf_t(x)\mathrm{d}\xi_t \\
 f_0(x) = & g(x) .
\end{split}
\end{equation}
for a process taking functions on $\mathbb{R}$ as values. Reading this equation as a famliy of ordinary SDE indexed by x, each one of which describes a real-valued process $f_t(x)$, one obtains the solution $f_t(x)=e^{x\xi_t-\frac{1}{2}x²t}g(x)$ (a family of geometric Brownian motions, indexed by x as diffusion constant). With regard to more complicated examples, the question arises whether such processes belong to certain function spaces, which would make them accessible methods of functional analysis. As for our example, $e^{x\xi_t-\frac{1}{2}x²t}=e^{-\frac{t}{2}(x-\frac{\xi_t}{t})² + \frac{\xi_t²}{2t}}$ is bounded as a function of x, so $g\in L^p(\mathbb{R})$ implies $f_t\in L^p(\mathbb{R})$ a.s., and also smoothness properties of  g are preserved. However, in the Hilbert space case $p=2$ we will be concerned with in the following, the computation
\begin{equation*}
 \mathbb{E}\|f_t\|^2 = \frac{1}{\sqrt{2\pi t}}\int e^{-\frac{y²}{2t}}\int f_t^2(x)\mathrm{d}x\mathrm{d}y
 = \frac{1}{\sqrt{2\pi t}}\int\int e^{-\frac{1}{2t}(y-2tx)²+tx²}g^2(x)\mathrm{d}y\mathrm{d}x
 =\int e^{tx²}g²(x)\mathrm{d}x
\end{equation*}
shows that $f_t$ remains in $L²(\Omega, L²(\mathbb{R}))$ only for very particular initial functions g (e.g. such with compact support). Otherwise, also the continuous dependence on the initial data cannot be expressed in terms of $\mathbb{E}\|\cdot\|²$. One can get rid of the troubles with the term $e^{tx²}$ by considering the process $\tilde f_t(x) := e^{-(1+c)\frac{t}{2}x²}f_t(x)$ with $c\ge0$. This modified process solves the equation
\begin{equation}\label{Diositeil}
 \begin{split}
  \mathrm{d}f_t(x) = & xf_t(x)\mathrm{d}\xi_t - (\frac{1}{2}+c)x²f_t(x)\mathrm{d}t\\
  f_0(x) = & g(x) .
 \end{split}
\end{equation}

The preceding observation suggests us to consider the equation
\begin{equation}\label{stochTeil}
\begin{split}
 \mathrm{d}\psi_t =& A\psi_t\mathrm{d}\xi_t - \left(\frac{1}{2}+c\right)A²\psi_t\mathrm{d}t \\
 \psi_0=& \psi\in\mathcal D(A²)
\end{split}
\end{equation}
with $c\ge0$ rather than $\mathrm{d}\psi_t = A\psi_t\mathrm{d}\xi_t$ also if we are dealing with more general self-adjoint operaors A on a Hilbert space H (the restriction to $\mathcal D(A²)$ will become clear). The formally expected solution
\begin{equation}\label{Loesstoch}
 \psi_t=e^{\xi_tA-(1+c)tA²}\psi
\end{equation}
is well-defined according to the functional calculus because $f(t,x,y):=e^{yx-(1+c)tx²}$, as a function of x with fixed values $t>0, y\in\mathbb{R}$ or $t=y=0$, is bounded on the spectrum of A. $\psi_t$ lies in $L²(\Omega,H)$ and depends continuously on an arbitrary given initial value $\psi$: For a proof let $\mu_\psi$ be the corresponding spectral measure. Then one computes
\begin{align*}
 \mathbb{E}\|\psi_t\|^2 =& \mathbb{E}\int f²(t,x,\xi_t)\mu_\psi(\mathrm{d}x)
 = \frac{1}{\sqrt{2\pi t}}\int e^{-\frac{y²}{2t}}\int f²(t,x,y)\mu_\psi(\mathrm{d}x)\mathrm{d}y \\
 =& \frac{1}{\sqrt{2\pi t}}\int\int e^{-\frac{1}{2t}(y-2tx)²-2ctx²} \mathrm{d}y\mu_\psi(\mathrm{d}x)
 = \int e^{-2ctx²} \mu_\psi(\mathrm{d}x) \\
 \le& \int 1 \mu_\psi(\mathrm{d}x) = \mu_\psi(\mathbb{R})=\|\psi\|²
\end{align*}
- note in particular that in the case $c=0$ even
\begin{equation}\label{isometrie}
\mathbb{E}\|\psi_t\|²=\|\psi\|²
\end{equation}
holds - and
\begin{equation}\label{stochAnf}
\begin{split}
 \mathbb{E}\|\psi_t-\psi\|² = \frac{1}{\sqrt{2\pi t}}\int e^{-\frac{y²}{2t}} \int (f(t,x,\xi_t)-1)²\mu_\psi(\mathrm{d}x)\mathrm{d}y \\
 = \int(1+e^{-2tx²}-2e^{-(1+c)tx²})\mu_\psi(\mathrm{d}x) \xrightarrow{t\to 0}0
\end{split}
\end{equation}
(since the integrand is majorized by a constant). By similar calculations on makes sure that, for $t>0$, $f(t,A,y)\psi$ is differentiable once w.r.t. t and twice w.r.t. y with the expected results. For example,
\begin{align*}
 &\left\| \frac{1}{h}(f(t+h,A,y)\psi-f(t,A,y)\psi) - [-(1+c)A²f(t,A,y)]\right\|^2 \\
 =& \int\frac{1}{h²}(f(t+h,x,y) - f(t,x,y) + h²(1+c)x²f(t,x,y))² \mu_\psi(\mathrm{d}x) \\
 =& \int e^{2yx-2(1+c)tx²}\left(\frac{e^{-(1+c)hx²}-1}{h}+(1+c)x²\right)^2 \mu_\psi(\mathrm{d}x) \xrightarrow{h\to0}0,
\end{align*}
because the estimate $|e^{-cx}-1|=c\int_0^x e^{-ct}\mathrm{d}t \le c\int_0^x\mathrm{d}t =cx$ valid for $c,x>0$ via
\begin{align*}
 & e^{2yx-2(1+c)tx²}\left(\frac{e^{-(1+c)hx²}-1}{h}+(1+c)x²\right)^2 \\
 \le& e^{2yx-2(1+c)tx²}\left(2\frac{(e^{-(1+c)hx²}-1)²}{h²}+2(1+c)²x^4\right)
 \le 4(1+c)²x^4e^{2yx-2(1+c)tx²}
\end{align*}
leads to a bounded and therefore $\mu_\psi$-integrable majorant.

According to the Itô formula, \eqref{Loesstoch} satisfies, at least for $t>0$ und $\epsilon\in(0,t)$, the equation
\begin{equation*}
 \psi_t = \psi_\epsilon+\int_\epsilon^tA\psi_s\mathrm{d}\xi_s - \int_\epsilon^t(\frac{1}{2}+c)A²\psi_s\mathrm{d}s
\end{equation*}
and it remains to be checked whether the terms on the right converge to their counterparts in \eqref{stochTeil} for $\epsilon\to0$. $\psi_\epsilon\to\psi$ has been shown in \eqref{stochAnf} and for $\psi\in \mathcal D(A²)$ we have, since in this case $e^{\xi_tA-(1+c)tA²}$ commutes with $A$ and $A²$,
\begin{align*}
 & \mathbb{E}\|\int_0^t(\frac{1}{2}+c)A²\psi_s\mathrm{d}s - \int_\epsilon^t(\frac{1}{2}+c)A²\psi_s\mathrm{d}s\|^2
 \le \epsilon\int_0^\epsilon\mathbb{E}\|(\frac{1}{2}+c)A²\psi_s\|^2\mathrm{d}s \\
 =& (\frac{1}{2}+c)^2\epsilon \int_0^\epsilon\mathbb{E}\|e^{\xi_sA-(1+c)sA²}A²\psi\|²\mathrm ds
 \le (\frac{1}{2}+c)^2\epsilon \int_0^\epsilon\mathbb{E}\|A²\psi\|²\mathrm ds \to0
\end{align*}
and
\begin{align*}
 \mathbb{E}\|\int_0^tA\psi_s\mathrm{d}\xi_s - \int_\epsilon^tA\psi_s\mathrm{d}\xi_s\|^2
 =\int_0^\epsilon\mathbb{E}\|A\psi_s\|²\mathrm ds
 \le \int_0^\epsilon\mathbb{E}\|A\psi\|²\mathrm ds \to0 .
\end{align*}

Let me emphasize once more that the conditions $c\ge0$ and $\psi\in\mathcal D(A²)$ which are required in the abstract ansatz can turn out to be far from optimal in concrete special cases - equation \eqref{Diositeil} is trivially solvable for arbirary c and even non-measurable initial values if one reads it, as performed for \eqref{geom}, as a family of ODE.

\section{General linear SDE with constant coefficients}

Let $\xi_t, \xi_t^1, \xi_t^2, \dots$ be independent Wiener processes on $(\Omega, \mathfrak{F}, \mathbb{P})$, $(\mathfrak{F}_t)_{t\ge0}$ the standard extension of the filtration generated by them, $\mathfrak{F}=\sigma(\mathfrak{F}_t,t\ge0)$, $\mathfrak{H}$ a separable $\mathbb{C}$-Hilbert space and $K, L_1, L_2, \dots$ operators on $\mathfrak{H}$ with dense common domain $\mathcal{D}$. We now develop the solution theory for the equation
\begin{equation}\label{linsdg}
\begin{split}
 \textnormal{d}\psi_t =& \sum_{j=1}^\infty L_j\psi_t\textnormal{d}\xi^j_t - K\psi_t\textnormal{d}t \\
 \psi_0 =& \psi\in\mathfrak{H},
\end{split}
\end{equation}
for an $\mathfrak{H}$-valued process $\psi_t(\omega)$ auf $\Omega$ with particular focus on the special case of the stochastic Schrödinger equation
\begin{equation}\label{stoschroedinger}
 \textnormal{d}\psi_t = (-iH-\frac{1}{2}A²)\psi_t\textnormal{d}t + A\psi_t\textnormal{d}\xi_t .
\end{equation}
Even in the case of bounded operators, no explicit solution formula for such equations seems to exist, but then at least the existence and uniqueness is ensured by global Lipschitz and growth conditions. In the general case, an approximation by bounded operators, as performed in \cite{Holevo}, at least leads to a weak solution. We present this construction and continue with the results of \cite{Mora} on the regularity of these solutions.

In the following we require the ``dissipativity condition''
\begin{equation}\label{Dissipativ}
 \sum_j \|L_j\psi\|² - 2\Re\langle K\psi,\psi\rangle \le c\|\psi\|² \textnormal{ für alle } \psi\in\mathcal{D} \textnormal{ und ein }c\ge0
\end{equation}
for the operators. It has the role of an a-priori growth estimate: Since a formal application of Itô's formula yields
\begin{equation*}
 \mathrm{d}e^{-ct}\|\psi_t\|² = 
 2e^{-ct}\sum\Re\langle\psi_t, L_j\psi_t\rangle\mathrm{d}\xi_t^j
 +e^{-ct} \underbrace{(- c\|\psi_t\|² - 2\Re\langle\psi_t,K\psi_t\rangle+\sum_j\|L_j\psi_t\|²)}_{\le0} \textnormal dt ,
\end{equation*}
it ensures that, if a sufficiently smooth solution $\psi_t$ of \eqref{linsdg} exists, $e^{-ct}\|\psi_t\|²$ is a supermaringale and, in particular, $\mathbb{E}\|\psi_t\|²\le e^{ct}\|\psi\|²$ gilt. Regularity is only analyzed for the special case $c=0$. Such equations, obviously comprising \eqref{stochTeil} from the last paragraph as well as the stochastic Schrödinger equation \eqref{stoschroedinger}, are called ``conservative''.

\subsection{Weak solutions}
In the sequel, by ``weak solution'' of \eqref{linsdg} we mean a process $\psi_t$ satisfying
\begin{equation*}
 \langle\phi,\psi_t\rangle = \langle\phi,\psi\rangle + \sum_j\int_0^t\langle L_j^*\phi,\psi_s\rangle\mathrm{d}\xi_s^j
 -\int_0^t\langle K^*\phi,\psi_s\rangle\mathrm{d}s
\end{equation*}
$\mathbb{P}$-a.s. for all $\phi$ in a dense set $\mathcal{D}^*$ (yet to be specified) on which $L_j^*$ and $K^*$ are defined. As for the existence of such a set, note that $K^*$ are in fact densely defined because $K_c:=K+\frac{c}{2}$ is accretive according to \eqref{Dissipativ} and thus closable, $\overline K=\overline{K_c}-\frac{c}{2}$ and $K^*\supset\overline K$. We assume the existence of a core $\mathcal{D}^*$ for $K^*$ such that all $L_j^*$ are defined on it and
\begin{equation}\label{2397}
 \forall \phi\in\mathcal{D}^*: \sum_{j=1}^\infty\|L_j^*\phi\|²<\infty .
\end{equation}

\begin{thm}\label{thm:Holevo}
 In the situation described above there exists a weak solution of \eqref{linsdg} with
 \begin{equation}\label{Wachstum}
  \mathbb{E}\|\psi_t\|²\le \|\psi\|²e^{ct} .
 \end{equation}
If $\overline{K+\frac{c}{2}id}$ is maximal accretive, the solution is unique under this condition.
\end{thm}
Basic facts about accretive operators can be found in \cite{Kato} or \cite{Reed2} and will be used without proof. Moreover, we need the following possibility of restricting the choice of test functions when dealing with weak convergence in $L²(\Omega)$:
\begin{lemm}\label{lemm:Chaos}
 Let $\xi_t^1, \dots, \xi_t^n$ be independent Wiener processes on $(\Omega, \mathfrak{F}, \mathbb{P})$, $(\mathfrak{F}_t)_{t\ge0}$ the standard extension of the filtration generated by them and $\mathfrak{F}=\sigma(\mathfrak{F}_t,t\ge0)$. We define recursively according to
\begin{equation*}
 I_0(t):=\chi_{[a,b]}(t) , \textnormal{ } I_{l+1}(t):=\sum_{j=1}^n\int_0^tI_l(s)a_j(s)\textnormal d\xi_s^j , \textnormal{ } I_l:=I_l(\infty)
\end{equation*}
stochastic processes, $a_j$ denoting arbitrary indicator functions of the form $\chi_{[a_j,b_j]}$. Then the random variables $I_0:=1$, $I_l:=I_l(\infty)$ (i.e. the integral in the recursive definition runs over the whole of $\mathbb{R}^+$) form a total subset of $L²(\Omega,\mathfrak F)$; the random variables $I_l(t)$ form a total subset of $L²(\Omega,\mathfrak F_t)$.
\end{lemm}

\begin{proof}
 According to a version of multidimensional Wiener chaos decomposition described in \cite{Shigekawa}, all iterated integrals of the form $\int_0^\infty\int_0^{t_{m-1}}\cdots\int_0^{t_1}\Phi(t_1,\dots,t_m)\textnormal d\xi_{t_1}^{i_1}\cdots\textnormal d\xi_{t_m}^{i_m}$ with $\Phi\in L^2(\mathbb{R}^m)$ form a total subset; since iterated integrals are continuous in the integrand, the choice of $\Phi$ can be restricted to the total subset of $L²(\mathbb{R}^m)$ constituted by products of the form $\chi_{[a_1,b_1]}(t_1)\cdots\chi_{[a_m,b_m]}(t_m)$. This proves the first part of the statement. In view of $I_l(t)=\mathbb{E}(I_l\mid\mathfrak F_t)$ and the $L²$-continuity of the projection $\mathbb E(\cdot\mid\mathfrak F_t)$, the $I_l(t)$ form a total subset of $L²(\Omega,\mathfrak F_t)$.
\end{proof}

\begin{proof}[Beweis von Satz \ref{thm:Holevo}]
 For technical simplicity we confine ourselves to $m\in\mathbb{N}$ instead of countably many Wiener processes.
We start by building modified Yosida approximations of $L_j$ and K by bounded operators: For this sake, let $\tilde K_c^*$ be an m-accretive extension of $K_c^*$ ($K_c^*=K^*+\frac{c}{2}$ is accretive because $K_c$ is) and $R_n:=(id+\frac{1}{n} \tilde K_c^*)^{-1}$ $(n\in\mathbb{N})$. Then $R_n$ is contractive, $s-\lim_{n\to\infty}R_n=id$ and $s-\lim_{n\to\infty}R_n^*=id$ (because, applying the general formula $B^*A^*\supset (AB)^*$ to the mutually inverse operators $R_n$ and $id+\frac{1}{n} \tilde K_c^*$, one obtains $R_n^*=(id+\frac{1}{n} \tilde K_c)^{-1}$, $\tilde K_c$ denoting an m-accretive extension of $K_c$). $L_j^n:=L_jR_n$ and $K_n:=R_n^*KR_n$ are bounded (they are closed and defined on the whole of $\mathfrak{H}$), unlike the usual Yosida approximation they satisfy the dissipativity condition ($\sum_j \|L_j^n\psi\|² - 2\Re\langle K^n\psi,\psi\rangle = \sum_j \|L_j(R_n\psi)\|² - 2\Re\langle K(R_n\psi),R_n\psi\rangle \le c\|R_n\psi\|² \le c\|\psi\|²$) and for all $\phi\in\mathcal{D}^*$ it holds true that $L_j^{n*}\phi\xrightarrow{n\to\infty}L_j^*\phi$ (because, for $\phi\in\mathcal D^*$, $L_j^{n*}\phi=R_n^*L_j^*\phi$ and $s-\lim R_n^*=id$) and $K^{n*}\phi\to K\phi$ (because $K^{n*}\phi=R_n^*K^*R_n\phi$, $K^*R_n$ isis the usual Yosida approximation of $K^*$, so it converges to $K^*$ on $\mathcal D^*$, then $R_n^*\to id$ and $\|R_n^*\|\le1$ imply the statement). As a consequence, the equations
 \begin{equation}\label{linsdgapp}
  \begin{split}
  \textnormal{d}\psi_t^n =& \sum_{j=1}^m L_j^n\psi_t^n\textnormal{d}\xi^j_t -  K^n\psi_t^n\textnormal{d}t \\
  \psi_0^n =& \psi
 \end{split}
 \end{equation}
 have a unique solution satisfying
 \begin{equation}\label{sfokm}
 \mathbb{E}\|\psi_t^n\|²\le e^{ct}\|\psi\|² .
 \end{equation}

This estimate shows that the $\psi_t^n$ are uniformly bounded as mappings from a time interval $[0,T]$ to $L²(\Omega,\mathfrak{H})$. We would like to show equicontinuity and apply the Arzela-Ascoli theorem. Since we have no clue concerning the limiting behaviour of expressions like $\mathbb{E}\|L_j^n\psi_t^n\|²$, we pass over to weak convergence in $L²(\Omega,\mathfrak{H})\cong L²(\Omega,\mathbb{C})\otimes\mathfrak{H}$, i.e. we choose $\xi\in L²(\Omega,\mathbb{C})$ and $\phi\in\mathcal{D}^*$ and consider the complex-valued functions $f^n(t):=\mathbb{E}\overline\xi\langle\phi,\psi_t^n\rangle$. Due to the boundedness of $\psi_t^n$,we only have to consider $\xi$ and $\phi$ from dense subsets of the unit speres, therefore we restrict to $\phi\in\mathcal{D}^*$; since $L²(\Omega,\mathbb C)$ and $\mathfrak H$ are separable (the first statement follows from the Wiener chaos decomposition), even an appropriate sequence $(\xi_l\phi_l)_{l\in\mathbb N}$ suffices. The $f^n(t)$ are uniformly bounded according to the Cauchy-Schwarz inequality and \eqref{sfokm} and equicontinuous according to
\begin{equation*}\label{gleichgradig1}
\begin{split}
 |f^n(t)-f^n(s)|² \le& \mathbb{E}|\xi|²\cdot\mathbb{E}|\langle\phi,\psi_t^n-\psi_s^n\rangle|²
 =\mathbb{E}|\sum_{j=1}^m\int_s^t\langle\phi,L_j^n\psi_\tau^n\rangle\mathrm{d}\xi_\tau^j - \int_s^t\langle\phi, K^n\psi_\tau^n\rangle\mathrm{d}\tau|² \\
 \le& 2\mathbb{E}|\sum_j\int_s^t\langle L_j^{n*}\phi,\psi_\tau^n\rangle\mathrm{d}\xi_\tau^j|²
 +2\mathbb{E}|\int_s^t\langle K^n\phi, \psi_\tau^n\rangle\mathrm{d}\tau|² \\
 \le& 2\sum_j\int_s^t\mathbb{E}|\langle L_j^{n*}\phi,\psi_\tau^n\rangle|²\mathrm{d}\tau
 +2|t-s|\int_s^t\mathbb{E}|\langle K^n\phi, \psi_\tau^n\rangle|²\mathrm{d}\tau \\
 \le& 2|t-s|\sum_j\|L_j^{n*}\phi\|²\sup_{\tau\in[s,t]}\mathbb{E}\|\psi_\tau\|²
 +2|t-s|²\|K^n\phi\|²\sup_{\tau\in[s,t]}\mathbb{E}\|\psi\|² \\
 \le& 2|t-s|\sum_j\|L_j^*\phi\|²e^{cT}\|\psi\|²
 +2|t-s|²\|K^n\phi\|²e^{cT}\|\psi_\tau\|² \\
 \le& C_{\phi,T}|t-s|\|\psi\|².
\end{split}
\end{equation*}
Therefore, there exists a subsequence $(n_k^1)_{k\in\mathbb N}$ of $(n)_{n\in\mathbb N}$ such that $\mathbb{E}\overline\xi_1\langle\phi_1,\psi_t^{n_k^1}\rangle$ converges uniformly in t, a subsequence $(n_k^2)$ of $(n_k^1)$ such that also $\mathbb{E}\overline\xi_2\langle\phi_2,\psi_t^{n_k^2}\rangle$ converges and so on. For the diagonal sequence $(n_k^k)=:(n_k)$ obviously $\mathbb{E}\overline\xi_j\langle\phi_j,\psi_t^{n_k^k}\rangle$ converges for all $j\in\mathbb N$ uniformly in t; since
\begin{equation*}
 L²(\Omega,\mathfrak{H})\ni\eta\mapsto\lim_{k\to\infty}\mathbb{E}\langle\eta,\psi_t^{n_k^k}\rangle\in\mathbb{R}
\end{equation*}
is a continuous functional, there exists a process $\psi_t$ with
\begin{equation*}
 \lim_{k\to\infty}\mathbb{E}\overline\xi\langle\phi,\psi_t^{n_k^k}\rangle = \mathbb{E}\overline\xi\langle\phi,\psi_t\rangle
\end{equation*}
for all $\xi,\phi$ uniformly in t. \eqref{sfokm} and the weak closedness of balls imply that this process satisfies \eqref{Wachstum}.

In order to show that this is the wanted solution, the convergence of the remaining summands in the ``weak form'' of \eqref{linsdgapp}, i.e. in
\begin{equation*}
 \mathbb{E}\overline\xi\langle\phi,\psi_t^{n_k}\rangle = \mathbb{E}\overline\xi\langle\phi,\psi\rangle
 + \mathbb{E}\overline\xi\langle\phi,\sum_{j=1}^m\int_0^tL_j^{n_k}\psi_s^{n_k}\mathrm{d}\xi_s^j\rangle
 - \mathbb{E}\overline\xi\langle\phi,\int_0^tK^{n_k}\psi_s^{n_k}\mathrm{d}s\rangle \textnormal{ resp.}
\end{equation*}
\begin{equation}\label{39}
 \mathbb{E}\overline\xi\langle\phi,\psi_t^{n_k}\rangle = \mathbb{E}\overline\xi\langle\phi,\psi\rangle
 + \sum_j\mathbb{E}\overline\xi\int_0^t\langle L_j^{n_k*}\phi,\psi_s^{n_k}\rangle\mathrm{d}\xi_s^j
 - \mathbb{E}\overline\xi\int_0^t\langle K^{n_k*}\phi,\psi_s^{n_k}\rangle\mathrm{d}s
\end{equation}
(strong Bochner integrals and their stochastic analogue allow for dragging bounded operators, in particular scalar products, inside),
to their counterparts in \eqref{linsdg} has to be checked.
To begin with,
\begin{align*}
 &|\mathbb{E}\overline\xi\int_0^t\langle K^{n_k*}\phi,\psi_s^{n_k}\rangle\mathrm{d}s
 -\mathbb{E}\overline\xi\int_0^t\langle K^*\phi,\psi_s\rangle\mathrm{d}s|
 = |\int_0^t\mathbb{E}\overline\xi\langle K^{n_k*}\phi,\psi_s^{n_k}\rangle\mathrm{d}s
 -\int_0^t\mathbb{E}\overline\xi\langle K^*\phi,\psi_s\rangle\mathrm{d}s| \\
 \le& |\int_0^t\mathbb{E}\overline\xi\langle K^{n_k*}\phi-K^*\phi,\psi_s^{n_k}\rangle\mathrm{d}s|
 +|\int_0^t\mathbb{E}\overline\xi\langle K^*\phi,\psi_s^{n_k}-\psi_s\rangle\mathrm{d}s| \\
 \le& (\mathbb{E}|\xi|^2)^{\frac{1}{2}} \|K^{n_k}\phi-K^*\phi\|\int_0^t\left(\mathbb{E}\|\psi_s^{n_k}\|^2\right)^{\frac{1}{2}} \mathrm ds +\hdots
 \to0;
\end{align*}
for the summands in the middle of \eqref{39} we restrict $\xi$ WLOG to the iterated stochastic integrals $I_l(t)$ specified in lemma \ref{lemm:Chaos}. This restriction to a total subset is justified by the fact that, according to
\begin{equation*}
 \mathbb{E}|\int_0^t\langle L_j^{n_k*}\phi,\psi_s^{n_k}\rangle\mathrm{d}\xi_s^j|^2
 =\int_0^t\mathbb{E}|\langle L_j^{n_k*}\phi,\psi_s^{n_k}\rangle|²\mathrm{d}s
 \le \|L_j^{n_k*}\phi\|^2 \int_0^t\mathbb{E}\|\psi_s^{n_k}\|²\mathrm{d}s
 \le CTe^{cT}\|\psi_0\|^2 ,
\end{equation*}
the sequences of stochastic integrals which are checked for convergence are bounded.
Then, similar as above and in addition using the Itô formula, the independence of the Wiener processes and the fact that, according to their recursive definition, all $I_l(t)$ are $L²$-martingales and have mean value 0, we conclude
\begin{align*}
 &|\mathbb{E}\overline{I_{l+1}(t)}\int_0^t\langle L_j^{n_k*}\phi,\psi_s^{n_k}\rangle\mathrm{d}\xi_s^j
 - \mathbb{E}\overline{I_{l+1}(t)}\int_0^t\langle L_j^*\phi,\psi_s\rangle\mathrm{d}\xi_s^j| \\
 =& |\sum_{h=1}^n\mathbb{E}\overline{\int_0^t a_h(s)I_l(s)\mathrm{d}\xi_s^h}\int_0^t \langle L_j^{n_k*}\phi,\psi_s^{n_k}\rangle\mathrm{d}\xi_s^j
 - \sum_h\mathbb{E}\overline{\int_0^t a_h(s)I_l(s)\mathrm{d}\xi_s^h}\int_0^t \langle L_j^*\phi,\psi_s\rangle\mathrm{d}\xi_s^j| \\
 =& |\mathbb{E}\int_0^t \overline{a_j(s)I_l(s)} \langle L_j^{n_k*}\phi,\psi_s^{n_k}\rangle\mathrm{d}s
 - \mathbb{E}\overline{\int_0^t a_j(s)I_l(s)} \langle L_j^*\phi,\psi_s\rangle\mathrm{d}s| \\
 \le& |\int_0^t\mathbb{E}\overline{a_j(s)I_l(s)} \langle L_j^{n_k*}\phi - L_j^*\phi,\psi_s^{n_k}\rangle\mathrm{d}s|
 + |\int_0^t\mathbb{E} \overline{a_j(s)I_l(s)} \langle L_j^*\phi,\psi_s^{n_k}-\psi_s\rangle\mathrm{d}s| \\
 \le& (\mathbb{E}|I_l(s)|^2)^\frac{1}{2}\|L_j^{n_k*}\phi - L_j^*\phi\|\int_0^t\mathbb{E}\|\psi_s^{n_k}\|^2\mathrm{d}s + \ldots \to0 .
\end{align*}

Now let $\psi_t^1$ and $\psi_t^2$ be two weak solutions of \eqref{linsdg} satisfying the same initial condition and \eqref{Wachstum} and define $\delta_t:=\psi_t^1-\psi_t^2$. For the unique solvability of \eqref{linsdg} it is sufficient to show $\mathbb{E}\overline{I_l(t)}\langle\phi,\delta_t\rangle=0$ for all $\phi\in\mathcal{D}^*$ and $I_l$ already specified. Since $\delta_0=0$,
\begin{equation}\label{733}
 \overline{I_0(t)}\langle\phi,\delta_t\rangle=\langle\phi,\delta_t\rangle
 = \sum_j\int_0^t\langle L_j^*\phi,\delta_s\rangle\mathrm{d}\xi_s^j
 -\int_0^t\langle K^*\phi,\delta_s\rangle\mathrm{d}s
\end{equation}
and, according to the product formula for Itô differentials,
\begin{equation}\label{734}
\begin{split}
 \overline{I_l(t)}\langle\phi,\delta_t\rangle =& \int_0^t \overline{I_{l-1}(s)}\langle\phi,\delta_s\rangle\sum_j \overline{a_j(s)}\textnormal d\xi_s^j
 +\int_0^t \overline{I_l(s)}\sum_j\langle L_j^*\phi,\delta_s\rangle \textnormal d\xi_s^j
 -\int_0^t \overline{I_l(s)}\langle K^*\phi,\delta_s\rangle \textnormal ds \\
 &+\int_0^t \overline{I_{l-1}(s)}\sum_j \overline{a_j(s)}\langle L_j^*\phi,\delta_s\rangle \textnormal{d}s
\end{split}
\end{equation}
for $l\ge1$. 
For the first time in this work, the integrands of the arising stochastic integrals contain products of two $L^2$-processes, so they are in $L^1$ according to Cauchy-Schwarz, but in general not in $L^2$. Therefore, the integrals are no longer defined in the usual $L^2$ framework via the Itô isometry, but only in the generalized sense described e.g. in \cite{Hackenbroch}, for which continuous paths are already sufficient. In general, such integrals are only local martingales. According to Holevo, the uniform boundedness of the $L^1$-norms of the integrands imply their martingale property. I have not found out why this is so - if it is true, then, taking mean values in \eqref{733} and \eqref{734} and defining $m_l(t):=\mathbb{E}\overline{I_l(t)}\delta_t$, one gets
\begin{equation}\label{735}
 \langle\phi,m_0(t)\rangle=-\int_0^t\langle K^*\phi,m_0(s)\rangle\textnormal ds
\end{equation}
and for $l\ge1$
\begin{equation}\label{736}
 \langle\phi,m_l(t)\rangle = \int_0^t\langle \sum_ja_j(s)L^{j*}\phi,m_{l-1}(s)\rangle \textnormal ds
 -\int_0^t\langle K^*\phi,m_l(s)\rangle \textnormal ds .
\end{equation}
Now let us prove inductively $m_l\equiv0$, which will of course imply $\mathbb{E}I_l(t)\langle\phi,\delta_t\rangle\equiv0$. Starting with $l=0$, we somehow have to bring the two terms $m_0$ in equation \eqref{735} together. For this sake, we take the Laplace transform, perform partial integration on the right-hand side und thus get rid of the additional integral in front of one $m_0$:
\begin{equation*}
 \int_0^te^{-\lambda s}\langle\phi,m_0(s)\rangle\textnormal ds
 = -\int_0^te^{-\lambda s}\int_0^s\langle K^*\phi,m_0(r)\rangle \textnormal dr\textnormal ds
 = -\frac{1}{\lambda}\int_0^te^{-\lambda s}\langle K^*\phi,m_0(s)\rangle \textnormal ds
\end{equation*}
Multiplying with $\lambda$ and dragging out the scalar products, we get
\begin{equation}\label{290}
 \langle\lambda\phi+K^*\phi,\int_0^te^{-\lambda s}m_0(s)\textnormal ds\rangle = 0 ,
\end{equation}
(the integrals still being understood as Bochner integrals).
Now let $\lambda>\frac{c}{2}$. According to the assumptions, $\overline{K+\frac{c}{2}}$ and thus also $K^*+\frac{c}{2}=\overline{K+\frac{c}{2}}^*$ is m-accretive, so $\Im(K^*+\lambda)=\mathfrak{H}$ and therefore $(K^*+\lambda)\mathcal D^*$ is dense in $\mathfrak H$. In view of that, \eqref{290} implies
\begin{equation*}
 \int_0^te^{-\lambda s}m_0(s)\textnormal ds=0
\end{equation*}
for $\lambda>\frac{c}{2}$. The Laplace transform is holomorphic, so this holds for all $\lambda>0$;
since a function is uniquely determined by its Laplace transform, $m_0\equiv0$ follows.
In the induction step $l-1\to l$, the first integral in \eqref{736} is 0 according to the induction hypothesis, so $m_l\equiv0$ follows from this equation along the same lines as $m_0\equiv0$.

Holevos uniqueness proof suffers from a (presumably only technical) gap. If, against expectation, it cannot be closed, let me point out that, at least under the additional assumptions in the next paragraph, the uniqueness of the strong solution of \eqref{stoschroedinger} is ascertained by its martingale property formulated in theorem \ref{thm:Martingal} (if there were two strong solutions $\psi_t^1$ and $\psi_t^2$ with the same initial value, $\|\psi_t^1-\psi_t^2\|^2$ was a martingale, so $\|\psi_0^1-\psi_0^2\|^2=0$ would imply $\mathbb{E}\|\psi_t^1-\psi_t^2\|^2=0$, so $\psi_t^1=\psi_t^2$ a.s.). Of course, the question whether there are additional weak solutions remains.
\end{proof}

\subsection{Regularität}
The weak solution $\psi_t$ is obviously strong if it does not leave the domains of $L_j$ and $K$ - or of another operator C with even smaller domain (then the operators in the weak form of the equation can be rolled back according to $\langle K^*\phi,\psi_t\rangle=\langle\phi,K\psi_t\rangle$ and the scalar products with test vectors $\phi$ be dropped). For a proof of regularity, it turns out to be appropriate to formulate conditions in terms of such a C rather than subjecting $L_j$ and $K$ to additional restrictions. ``C-strong solution'' of a conservative equation of the form \eqref{linsdg} then denotes an adapted, a.s. D(C)-valued process with continuous paths that satisfies the equation and, moreover, the estimates $\mathbb{E}\|\psi_t\|²\le \mathbb{E}\|\psi\|²$ and $\sup_{0\le s\le t}\mathbb{E}\|C\psi_t\|²<\infty$ for all $t\ge0$.

\begin{thm}\label{starkeLoes}
 Consider a conservative equation of the form \eqref{linsdg} with $D(K)\subset D(L_j)$ - that is to say, the operators have to satisfy
 \begin{equation}\label{konservativ}
 \sum_j \|L_j\psi\|² - 2\Re\langle K\psi,\psi\rangle =0
 \end{equation}
 for all $\psi\in D(K)$.
 Let C be a positive definite self-adjoint ``reference operator'' with the following properties:
 \begin{itemize}
  \item $D(C)\subset D(K)\cap D(K^*)$.
  \item There exists an orthonormal basis $e_1,e_2,\dots$ of  $\mathfrak{H}$ consisting of elements from $D(C)$ such that $\forall n\in\mathbb{N}: \sum_j \|L_j^*e_n\|²<\infty$.
  \item Let $P_n$ be the orthogonal projection onto $\mathfrak{H}_n:=\langle e_1,\dots,e_n\rangle$. Then there are constants $\alpha,\beta\ge 0$ such that
  \begin{equation}
   2\Re\langle Cx,CP_nKx\rangle + \sum_j\|CP_nL_jx\|² \le \alpha\left(\|Cx\|²+\|x\|²+\beta\right)
  \end{equation}
  holds for all $n\in\mathbb{N}$ and $x\in\mathfrak{H}_n$.
  \item $\sup_{n\in\mathbb{N}}\|CP_nx\|\le \|Cx\|$ for all $x\in D(C)$.
 \end{itemize}
 Finally, let $\psi$ be a $\mathfrak{F}_0$-measurable, a.s. $D(C)$-valued random variable with $\mathbb{E}\|\psi\|_C^2<\infty$ (as usual, we write $\|x\|_C:=\sqrt{\|x\|²+\|Cx\|²}$).
 Then there exists a unique C-strong solution of \eqref{linsdg} which satisfies the above mentioned estimates for all $t\ge0$. This solution also satisfies
 \begin{equation*}
  \mathbb{E}\|C\psi_t\|² \le e^{\alpha t}\left(\mathbb{E}\|C\psi\|²+\alpha t(\mathbb{E}\|\psi\|²+\beta)\right) .
 \end{equation*}

\end{thm}

\begin{thm}\label{thm:Martingal}
 With $\psi_t$ denoting the solution from theorem \ref{starkeLoes}, $\|\psi_t\|²$ is a martingale, in particular $\mathbb{E}\|\psi_t\|²=\|\psi\|²$.
\end{thm}

\subsection{Comment on the form of the theorems}
The conditions in the theorems appear to be formulated ad hoc. The authors comment that they take inspirations from ``quantum stochastic calculus'' - however, without further explanation, the conditions still look like a list of possible problems which one might encounter when dealing with these equations, but the real essence of which has not been inquired yet. However, the existence of  weak solutions is ensured for many concrete applications - e.g. if $\mathfrak H=L²(\mathbb R^n)$ and the operators $K$ and $L_j$ are ``built out ouf'' differential and multiplication operators, then ``usually'' $\mathcal D=\mathcal D^*=C_c^\infty$ will be suitable domains. Finding a reference operator C seems to be a more touchy business, which is not solvable by a general algorithm. In \cite{Mora} Mora and Rebolledo have found such an operator i.a. for the stochastic Schrödinger equation
\begin{equation*}
 \textnormal d\psi_t = -\frac{i}{2}\triangle\psi_t + x\psi_t\textnormal d\xi_t
 - \frac{1}{2}x^2\psi_t\textnormal dt ,
\end{equation*}
tracing back to Diosi (\cite{Diosi2}), namely an arbitrary natural power of $N:=x^2-\frac{1}{2}\triangle-id$, an operator which is self-adjoint on
\begin{equation}\label{Domref}
\{f\in H^2(\mathbb{R})|x²f\in L²(\mathbb{R})\} .
\end{equation}
Even if this result looks satisfying, one has to mention the articles \cite{Kolokoltsov, Kolb}, in which the authors have tailored a special method for this equation (transformation to the equation of the non-self adjoint harmonic oscillator) and found an integral kernel which provides an explicit strong solution for arbitrary bounded $L^2$ initial conditions. So one might guess that our restrictions are due to the particular ansatz rather than the actual behaviour of the equations.

\chapter{Product formula}\label{ch:Produkt}

\section{Linear SDE in $\mathbb{R}^m$}

Define a Wiener process $\xi_t$ on a probability space $(\Omega,\mathfrak{A},\mathbb{P})$ and consider the m-dimensional SDE
\begin{equation}\label{Gleichung}
\begin{split}
 \textnormal{d}X_t=&AX_t\textnormal{d}t+BX_t\textnormal{d}\xi_t \\
 X_0=&X
\end{split}
\end{equation}
with $A,B\in \mathbb{R}^{m\times m}$ and $X\in \mathbb{R}^m$. The deterministic Trotter formula and the above mentioned result of \cite{Gough} suggest that its solution flow $C_{s,t}$ can be obtained from the equations
\begin{equation}\label{det}
 \textnormal{d}X_t=AX_t\textnormal{d}t
\end{equation}
and
\begin{equation}\label{stoch}
 \textnormal{d}X_t=BX_t\textnormal{d}\xi_t ,
\end{equation}
resp. their solution flows
\begin{equation}\label{det'}
 A_{s,t}=e^{(t-s)A}
\end{equation}
and
\begin{equation}\label{stoch'}
 B_{s,t}=e^{(\xi_t-\xi_s)B-\frac{(t-s)²}{2}B²}
\end{equation}
via the product formula
\begin{equation}\label{Produkt}
 C_{0,t}=\lim_{n\rightarrow\infty}(\prod_{k=0}^{n-1}e^{\frac{t}{n}A}B_{\frac{k}{n}t,\frac{k+1}{n}t}) .
\end{equation}
Here and in the following, products of operators are to be understood in the order of their action on a vector, i.e.
\begin{equation*}
 \prod_{k=1}^nH_k:=H_n\cdots H_1 .
\end{equation*}

In case that single random variables $C_{0,t}$ can be approximated in such a way, in contrast to the deterministic case the additional question arises whether the right-hand side in \eqref{Produkt} is also a good approximation of $(C_{0, t})_{t\ge 0}$ when considered as a stochastic process, i.e. one should also compare their autocorrelations or inquire the convergence of whole paths. For this purpose, it turned out to be useful to restrict to finite time intervals [0,T] and to replace the process $(\prod_{k=0}^{n-1}e^{\frac{t}{n}A}B_{\frac{k}{n}t,\frac{k+1}{n}t})_{t\ge 0}$, which makes the same number of approximation steps for all times, either by a process $(f_{n,T}(t))_{t\in [0,T]}$, which is defined as $f_{n,T}\left(\frac{k}{n}T\right):=\prod_{j=0}^{k-1}e^{\frac{t}{n}A}B_{\frac{j}{n}T,\frac{j+1}{n}T}$ on lattice points $t\in \{0, \frac{T}{n}, \dots, T\}$ and constantly extended in between, that is to say $f_{n,T}\left(\frac{k}{n}T+\tau\right):=f_{n,T}\left(\frac{k}{n}T\right)$ for $t=\frac{k}{n}T+\tau$ with $\tau\in(0,\frac{T}{n})$, or by a process $g_{n,T}$ which is defined in the same way on the lattice points and interpolated continuously according to $g_{n,T}\left(\frac{k}{n}T+\tau\right):=e^{\tau A}B_{\frac{k}{n}T,\frac{k}{n}T+\tau}\prod_{j=0}^{k-1}e^{\frac{t}{n}A}B_{\frac{j}{n}T,\frac{j+1}{n}T}$. In order to avoid long case differentiations, one can write this definition by use of the Gauß bracket as
\begin{equation}
 g_{n,T}(t):=e^{\left(t- \lfloor\frac{nt}{T}\rfloor\frac{T}{n}\right)A}B_{\lfloor\frac{nt}{T}\rfloor\frac{T}{n}, t} \prod_{k=0}^{\lfloor\frac{nt}{T}\rfloor-1}e^{\frac{t}{n}A}B_{\frac{k}{n}T,\frac{k+1}{n}T} .
\end{equation}

As will be shown, the paths of these processes converge to the ones of $C_{0,\cdot}$ uniformly in the quadratic mean. Thus, or notion of convergence is slightly improved as  compared with \cite{Gough}, where uniformity is missing. Moreover, we have noted that, on the analogy of Chernoff's theorem and also beyond the results of \cite{Gough}, the factors $F(s,t):=e^{(t-s)A}B_{s,t}$ appearing in \eqref{Produkt} can be altered, as long as they keep solving \eqref{Gleichung} ``for infinitesimal time intervals'' when regarded as stochastic processes in the time parameter t with fixed s:

\begin{thm}\label{Satz1}
 Let $F: \{(x,y)\in \mathbb{R}²:y\ge x\ge 0\} \rightarrow \textnormal{L}^2(\Omega, \mathbb{R}^{m\times m})$ be a mapping with the following properties: \\
 (i)\ \ $F(s,t)=1+(t-s)A+(\xi_t-\xi_s)B+o(t-s)$ or \\
 (i)'\ $F(s,t)=1+\int_s^t AF(s,\tau)\textnormal d\tau+\int_s^t BF(s,\tau)\textnormal d\xi_\tau+o(t-s)$
 für $t\rightarrow s$ \\
 (ii) For $\Delta>0$ and $u\ge s+\Delta$, $F(s,s+\Delta)$ and $F(u,u+\Delta)$ are i.i.d. \\
 Define $f_{n,T}(t):=\prod_{k=0}^{\lfloor\frac{nt}{T}\rfloor-1}F(\frac{k}{n}T,\frac{k+1}{n}T)$ for $t\in [0,T]$, i.e.
 $f_{n,T}(t)=\prod_{j=0}^{k-1}F(\frac{j}{n}T,\frac{j+1}{n}T)$ if $t\in \left[\frac{k}{n}T,\frac{k+1}{n}T\right[$, and
 $g_{n,T}:=F\left(\lfloor\frac{nt}{T}\rfloor\frac{T}{n},t\right)\prod_{k=0}^{\lfloor\frac{nt}{T}\rfloor-1}F(\frac{k}{n}T,\frac{k+1}{n}T)$. \\
 a) (i) and (ii) imply $\mathbb{E}\sup_{t\in [0,T]}\|f_{n,T}(t)-C_{0,t}\|²\xrightarrow{n\to\infty}0$. \\
 b) (i)' and (ii) imply $\mathbb{E}\sup_{t\in [0,T]}\|g_{n,T}(t)-C_{0,t}\|²\xrightarrow{n\to\infty}0$.
\end{thm}

\begin{rem}
 (i) and (i') are not special cases of one another because $\int_s^t AF(s,\tau)\textnormal d\tau=(t-s)A+o(t-s)$, but in general only $\int_s^t BF(s,\tau)\textnormal d\xi_\tau=(\xi_t-\xi_s)B+o(\sqrt{t-s})$ holds.
\end{rem}
\begin{rem}\label{rem:Vorfaktor}
 Under the additional assumption that $\mathbb{E}\sup_{t\in[0,\frac{T}{n}]}\|o(t)\|^2\to0$ for $n\to\infty$, the convergence of $f_{n,T}$ continues to hold true under condition (i)' and $g_{n,T}$ also converges under (i): Assume that the convergence formulated in the theorem and, in particular, the boundedness in the form $\sup_{n\in\mathbb{N}}\mathbb{E}\sup_{0\le t\le T}\|f_{n,T}(t)\|^2<\infty$ and the same for $g_{n,T}$ are already proven. Then, since
 \begin{align*}
  \mathbb{E}\sup_{t\in[0,T]}\|g_{n,T}(t)-f_{n,T}(t)\|^2
 =&\mathbb{E}\sup_t\|\left(F\left(\left\lfloor\frac{nt}{T}\right\rfloor\frac{T}{n},t\right)-1\right) f_{n,T}(t)\|^2 \\
 \le& \mathbb{E}\sup_t\|\left(F\left(\left\lfloor\frac{nt}{T}\right\rfloor\frac{T}{n},t\right)-1\right)\|^2
 \mathbb{E}\sup_t\|f_{n,T}(t)\|^2 ,
 \end{align*}
 one would only have to show
$\lim_{n\to\infty}\mathbb{E}\sup_t\|\left(F\left(\left\lfloor\frac{nt}{T}\right\rfloor\frac{T}{n},t\right)-1\right)\|^2 =0$
which, considering e.g. for case (i)
\begin{align*}
 &\mathbb{E}\sup_t\|F\left(\left\lfloor\frac{nt}{T}\right\rfloor\frac{T}{n},t\right)-1\|^2
 =\mathbb{E}\sup_t\|(t-\left\lfloor\frac{nt}{T}\right\rfloor\frac{T}{n})A 
 +(\xi_t-\xi_{\left\lfloor\frac{nt}{T}\right\rfloor\frac{T}{n}})B + o(t-\left\lfloor\frac{nt}{T}\right\rfloor\frac{T}{n})\|^2 \\
 \le& 3\frac{T²}{n²}\|A\|² + 3\mathbb{E}\sup_{t}(\xi_t-\xi_{\left\lfloor\frac{nt}{T}\right\rfloor\frac{T}{n}})^2\|B\|^2  + o(1) ,
\end{align*}
boils down to $\mathbb{E}\sup_{t}(\xi_t-\xi_{\left\lfloor\frac{nt}{T}\right\rfloor\frac{T}{n}})^2 \to0$ and is shown in the following proof from \eqref{f73} on.
\end{rem}

\begin{proof}[Beweis von Satz \ref{Satz1}]
In order to check the plausibility of the theorem and  as a technical tool we first show that $f_{n,T}$ and $g_{n,T}$ are bounded in $\textnormal{L}^2(\Omega, \mathbb{R}^{m\times m})$ on $[0,T]$ by constants not depending on n. A crude application of the triangle inequality for the $\textnormal{L}^2(\Omega, \mathbb{R}^{m\times m})$-norm of the single factors along with (ii) leads to
\begin{equation*}\label{falsch}
\begin{split}
 \sqrt{\mathbb{E}\|f_{n,T}(t)\|²} \le
 &\left(\mathbb{E}\prod_{k=0}^{\lfloor\frac{nt}{T}\rfloor-1}\|F(\frac{k}{n}T,\frac{k+1}{n}T)\|²\right)^{\frac{1}{2}} =
 \prod_{k=0}^{\lfloor\frac{nt}{T}\rfloor-1}\left(\mathbb{E}\|F(\frac{k}{n}T,\frac{k+1}{n}T)\|²\right)^\frac{1}{2} = \\
 = &\sqrt{\left(\mathbb{E}\|F(0,\frac{T}{n})\|²\right)}^{\lfloor\frac{nt}{T}\rfloor} \le
 \left(1 + \frac{T}{n}\|A\| + \sqrt{\frac{T}{n}}\|B\| + \epsilon\frac{T}{n}\right)^n
\end{split}
\end{equation*}
and, since $\left(1+\frac{c}{\sqrt{n}}\right)^n \approx e^{c\sqrt{n}}$, does not ensure boundedness, but this can be improved by reading the norm on $\mathbb{R}^{m\times m}$ as the Frobenius norm generated by $\langle A,B\rangle:=tr(A^tB)$ (other possible norms would be equivalent anyway) and replacing the triangle inequality by an exact computation. In case a) this leads to
\begin{equation}\label{beschr}
\begin{split}
 &\mathbb{E}\|f_{n,T}(t)\|² \le \mathbb{E}\prod_{k=0}^{\lfloor\frac{nt}{T}\rfloor-1}\|F(\frac{k}{n}T,\frac{k+1}{n}T)\|² =
 \prod_{k=0}^{\lfloor\frac{nt}{T}\rfloor-1}\mathbb{E}\|F(\frac{k}{n}T,\frac{k+1}{n}T)\|² \\
 =& \prod_{k=0}^{\lfloor\frac{nt}{T}\rfloor-1} \left(\|1+\frac{T}{n}A\|²+\mathbb{E}\|(\xi_{\frac{k+1}{n}T}-\xi_{\frac{k}{n}T})B\|²+
 \mathbb{E}\|o_k\left(\frac{1}{n}\right)\|²+ \color{white}\right) \\
 &\color{white}\left( \color{black} +2\mathbb{E}\langle 1+\frac{T}{n}A,(\xi_{\frac{k+1}{n}T}-\xi_{\frac{k}{n}T})B\rangle+
 2\mathbb{E}\langle (\xi_{\frac{k+1}{n}T}-\xi_{\frac{k}{n}T})B,o_k\left(\frac{1}{n}\right)\rangle+
 2\mathbb{E}\langle 1+\frac{T}{n}A,o_k\left(\frac{1}{n}\right)\rangle\right) \\
 \le& \prod_{k=0}^{\lfloor\frac{nt}{T}\rfloor-1} \left(\left(1+\frac{T}{n}\|A\|\right)²+\mathbb{E}(\xi_{\frac{k+1}{n}T}-\xi_{\frac{k}{n}T})²\|B\|²+
 o\left(\frac{1}{n²}\right)
 +2\mathbb{E}(\xi_{\frac{k+1}{n}T}-\xi_{\frac{k}{n}T})\langle 1+\frac{T}{n}A,B\rangle  \color{white}\right) \\
 &\color{white}\left( \color{black} +2\left(\mathbb{E}\|(\xi_{\frac{k+1}{n}T}-\xi_{\frac{k}{n}T})B\|²\right)^{\frac{1}{2}}
 \left(\mathbb{E}\|o_k\left(\frac{1}{n}\right)\|²\right)^{\frac{1}{2}}+
 2\|1+\frac{T}{n}A\|\left(\mathbb{E}\|o_k\left(\frac{1}{n}\right)\|²\right)^{\frac{1}{2}}\right) \\
 =& \left(\left(1+\frac{T}{n}\|A\|\right)²+\frac{T}{n}\|B\|²+
 o\left(\frac{1}{n²}\right)+0+
 2\sqrt{\frac{T}{n}}\|B\|
 o\left(\frac{1}{n}\right)+
 2\|1+\frac{T}{n}A\|o\left(\frac{1}{n}\right)\right)^{\lfloor\frac{nt}{T}\rfloor-1} \\
 \le& \left(1+\mathcal{O}\left(\frac{1}{n}\right)\right)^n .
\end{split}
\end{equation}
This term remains bounded because $\left(1+\frac{c}{n}\right)^n\to e^c$. Now we have avoided a $\frac{1}{\sqrt{n}}$ summand within the bracket because $\mathbb{E}\langle 1, (\xi_{\frac{k+1}{n}T}-\xi_{\frac{k}{n}T})B\rangle = 0$ was computed instead of being estimated by Cauchy-Schwarz. The same trick works in case b), $(\xi_{\frac{k+1}{n}T}-\xi_{\frac{k}{n}T})B$ now being replaced by $\int_{\frac{k}{n}T}^{\frac{k+1}{n}T}BF\left(\frac{k}{n}T, \tau\right)\mathrm{d}\xi_\tau$; as for the other terms, it suffices to estimate the integrals via $\int_s^t a_\tau\textnormal{d}\tau=(t-s)a_s+o(t-s)$ and $\int_s^t a_\tau\textnormal{d}\xi_\tau=(\xi_t-\xi_s)a_0+o(\sqrt{t-s})$ (valid in $L²$) and products of them via Cauchy-Schwarz.

The observation that $f_{n,T}$ and $g_{n,T}$ deviate only a little from the flow $C_{0,t}$ can be formulated in terms of integral equations similar to
\begin{equation*}
 C_{0,t}=1+\int_0^t AC_{0,s}\textnormal ds + \int_0^t BC_{0,s}\textnormal{d}\xi_s :
\end{equation*}
In Fall a), one gets
 \begin{equation}\label{Teleskop}
 \begin{split}
  f_{n,T}(t) &= f_{n,T}\left(\left\lfloor\frac{nt}{T}\right\rfloor\frac{T}{n}\right) =
  1+\sum_{k=0}^{\lfloor\frac{nt}{T}\rfloor-1}\left(f_{n,T}\left(\frac{k+1}{n}T\right)-f_{n,T}\left(\frac{k}{n}T\right)\right) \\
  &=1+\sum_{k=0}^{\lfloor\frac{nt}{T}\rfloor-1}\left(F\left(\frac{k}{n}T,\frac{k+1}{n}T\right)-1\right)f_{n,T}\left(\frac{k}{n}T\right)  \\
  &=1+\sum_{k=0}^{\lfloor\frac{nt}{T}\rfloor-1}\left(\frac{T}{n}A+(\xi_{\frac{k+1}{n}T}-\xi_{\frac{k}{n}T})B+
  o_k\left(\frac{1}{n}\right)\right)f_{n,T}\left(\frac{k}{n}T\right) \\
  &=1+\sum_{k=0}^{\lfloor\frac{nt}{T}\rfloor-1}\left(\int_{\frac{k}{n}T}^{\frac{k+1}{n}T}Af_{n,T}\left(\frac{k}{n}T\right)\textnormal{d}s
  +\int_{\frac{k}{n}T}^{\frac{k+1}{n}T}Bf_{n,T}\left(\frac{k}{n}T\right)\textnormal{d}\xi_s
  +o_k\left(\frac{1}{n}\right)f_{n,T}\left(\frac{k}{n}T\right)\right) \\
  &=1+\int_0^{\lfloor\frac{nt}{T}\rfloor\frac{T}{n}}Af_{n,T}(s)\textnormal{d}s
  +\int_0^{\lfloor\frac{nt}{T}\rfloor\frac{T}{n}}Bf_{n,T}(s)\textnormal{d}\xi_s
  +\sum_{k=0}^{\lfloor\frac{nt}{T}\rfloor-1}o_k\left(\frac{1}{n}\right)f_{n,T}\left(\frac{k}{n}T\right) .
 \end{split}
 \end{equation}
Using this representation, we aim at comparing $f_{n,T}$ and $C_{0,\cdot}$ by means of the Gronwall inequality. For this end, consider $d_{n,T}(S):=\mathbb{E}\sup_{0\le t\le S}\|f_{n,T}(t)-C_{0,t}\|²$ for arbitrary $S\le T$. In
\begin{align*}
 d_{n,T}(S) \le&
 5\mathbb{E}\sup_{0\le t\le S}\left\|\int_0^{\lfloor\frac{nt}{T}\rfloor\frac{T}{n}}A(f_{n,T}(s)-C_{0,s})\textnormal{d}s\right\|^2
  +5\mathbb{E}\sup_{0\le t\le S}\left\|\int_0^{\lfloor\frac{nt}{T}\rfloor\frac{T}{n}} B(f_{n,T}(s)-C_{0,s})\textnormal{d}\xi_s\right\|^2 \\
  &+5\mathbb{E}\sup_{0\le t\le S}\left\|\int_{\lfloor\frac{nt}{T}\rfloor\frac{T}{n}}^t AC_{0,s}\textnormal{d}s\right\|^2
  +5\mathbb{E}\sup_{0\le t\le S}\left\|\int_{\lfloor\frac{nt}{T}\rfloor\frac{T}{n}}^t BC_{0,s}\textnormal{d}\xi_s\right\|^2 \\
  &+5\mathbb{E}\sup_{0\le t\le S}\left\|\sum_{k=0}^{\lfloor\frac{nt}{T}\rfloor-1} o_k\left(\frac{1}{n}\right)f_{n,T}\left(\frac{k}{n}T\right)\right\|^2 ,
\end{align*}
the first summand can be estimated according to
\begin{align*}
 &\mathbb{E}\sup_{0\le t\le S}\left\|\int_0^{\lfloor\frac{nt}{T}\rfloor\frac{T}{n}}A(f_{n,T}(s)-C_{0,s})\textnormal{d}s\right\|² \le
 \mathbb{E}\sup_{0\le t\le S}\lfloor\frac{nt}{T}\rfloor\frac{T}{n} \int_0^{\left\lfloor\frac{nt}{T}\right\rfloor\frac{T}{n}}\|A(f_{n,T}(s)-C_{0,s})\|²\textnormal{d}s \\
 \le& T\|A\|²\mathbb{E}\int_0^{S}\|f_{n,T}(s)-C_{0,s}\|²\textnormal{d}s =
 T\|A\|²\int_0^{S}\mathbb{E}\|f_{n,T}(s)-C_{0,s}\|²\textnormal{d}s \le
 T\|A\|²\int_0^{S}d_{n,T}(s)\textnormal{d}s ,
\end{align*}
for the ensuing stochastic integral one additionally uses Doob's inequality and Itô's formula and gets
\begin{gather*}
 \mathbb{E}\sup_{0\le t\le S}\left\|\int_0^{\lfloor\frac{nt}{T}\rfloor\frac{T}{n}} B(f_{n,T}(s)-C_{0,s})\textnormal{d}\xi_s\right\|² \le
 4\mathbb{E}\left\|\int_0^{\lfloor\frac{nS}{T}\rfloor\frac{T}{n}} B(f_{n,T}(s)-C_{0,s})\textnormal{d}\xi_s\right\|² = \\
 =4\int_0^{\lfloor\frac{nS}{T}\rfloor\frac{T}{n}}\mathbb{E}\|B(f_{n,T}(s)-C_{0,s})\|²\textnormal{d}s \le
 4\|B\|²\int_0^S d_{n,T}(s)\textnormal{d}s .
\end{gather*}
In the remaining summands $d_{n,T}$ does not appear - if they converged to 0 uniformly in S for $n\to\infty$, then
\begin{equation*}
 d_{n,T}(S) \le (T\|A\|²+4\|B\|²)\int_0^Sd_{n,T}(s)\textnormal{d}s + o(1)
\end{equation*}
would imply at first $d_{n,T}(T)\le o(1)e^{(T\|A\|²+4\|B\|²)T}$ and for $n\to\infty$ the assertion of the theorem.

Whereas

\begin{gather*}
 \mathbb{E}\sup_{0\le t\le S}\left\|\int_{\lfloor\frac{nt}{T}\rfloor\frac{T}{n}}^t AC_{0,s}\textnormal{d}s\right\|² \le
 \frac{T}{n}\|A\|²\int_0^T\mathbb{E}\|C_{0,s}\|²\textnormal{d}s =
 \mathcal{O}\left(\frac{1}{n}\right)
\end{gather*}
is straightforward and also the ``little o'' term, given that $o_k\left(\frac{1}{n}\right)$ is independent of $f_{n,T}\left(\frac{k}{n}T\right)$ and $\mathbb{E}\|o_k\left(\frac{1}{n}\right)\|²=o(\frac{1}{n²})$ holds, satisfies
\begin{equation}\label{Rest}
 \begin{split}
 & \mathbb{E}\sup_{0\le t\le S} \left\|\sum_{k=0}^{\lfloor\frac{nt}{T}\rfloor-1}o_k\left(\frac{1}{n}\right)f_{n,T}\left(\frac{k}{n}T\right)\right\|² \le
 n\mathbb{E}\sup_{0\le t\le T}\sum_{k=0}^{\lfloor\frac{nt}{T}\rfloor-1}\left\|o_k\left(\frac{1}{n}\right)f_{n,T}\left(\frac{k}{n}T\right)\right\|² \\
 \le & n\sum_{k=0}^{n-1}\mathbb{E}\left\|o_k\left(\frac{1}{n}\right)f_{n,T}\left(\frac{k}{n}T\right)\right\|² \le
 no\left(\frac{1}{n²}\right)\sum_{k=0}^{n-1}\mathbb{E}\left\|f_{n,T}\left(\frac{k}{n}T\right)\right\|² \\
 \le & o(1)\sup_{t\in [0,T], n\in\mathbb{N}}\mathbb{E}\|f_{n,T}(t)\|²=o(1) ,
 \end{split}
\end{equation}
the remaining proof of $\lim_{n\to\infty}\mathbb{E}\sup_{0\le t\le S} \left\| \int_{\lfloor\frac{nt}{T}\rfloor\frac{T}{n}}^t BC_{0,s}\textnormal{d}\xi_s\right\|²=0$ requires somewhat more effort. It holds that
\begin{equation}\label{f73}
\begin{split}
 &\mathbb{E}\sup_{0\le t\le S}\left\| \int_{\lfloor\frac{nt}{T}\rfloor\frac{T}{n}}^t BC_{0,s}\textnormal{d}\xi_s\right\|²
 \le\mathbb{E}\sup_{0\le t\le T} \left\| \int_{\lfloor\frac{nt}{T}\rfloor\frac{T}{n}}^t BC_{0,s}\textnormal{d}\xi_s\right\|²\\
 =&\mathbb{E}\sup_{0\le k\le n-1 , 0\le t<\frac{1}{n}} \left\|M_{\frac{k}{n}T+t}-M_{\frac{k}{n}T}\right\|²
\end{split}
\end{equation}
with the $L²$-martingale $M_t:=\int_0^t BC_{0,s}\textnormal{d}\xi_s$. First let n be fixed. By the Doob inequality, all increments $M_{\frac{k}{n}T+t}-M_{\frac{k}{n}T}$ with fixed k and $t\in [0,\frac{1}{n}]$ can be bounded by moments of $M_{\frac{k+1}{n}T}-M_{\frac{k}{n}T}$ and increments over disjoint time intervals are independent. If $M_t$ were replaced by the (one-dimensional) Wiener process $\xi_t$, they would have identical normal distributions and all this would give rise to an explicit calculation:
\begin{gather*}
 \mathbb{E}\sup_{0\le k\le n-1, 0\le t<\frac{1}{n}} \left|\xi_{\frac{k}{n}T+t}-\xi_{\frac{k}{n}T}\right|²=
 \int_0^\infty \mathbb{P}\left(\sup_{0\le k\le n-1 , 0\le t<\frac{1}{n}} \left|\xi_{\frac{k}{n}T+t}-\xi_{\frac{k}{n}T}\right|²>s\right)\textnormal{d}s= \\
 =\int_0^\infty \left(1-\mathbb{P}\left(\forall k\in \{0,\dots, n-1\}: \sup_{0\le t<\frac{1}{n}} \left|\xi_{\frac{k}{n}T+t}-\xi_{\frac{k}{n}T}\right|²\le s\right)\right)\textnormal{d}s = \\
 =\int_0^\infty \left(1-\mathbb{P}\left(\sup_{0\le t<\frac{1}{n}} \left|\xi_t-\xi_0\right|²\le s\right)^n\right)\textnormal{d}s =
 \int_0^\infty \left(1-\left(1-\mathbb{P}\left(\sup_{0\le t<\frac{1}{n}} |\xi_t|²> s\right)\right)^n\right)\textnormal{d}s \le \\
 \le \int_0^\infty \left(1-1\wedge\left|1-\frac{\mathbb{E}\xi_\frac{1}{n}^4}{s²}\right|^n\right)\textnormal{d}s =
 \int_0^\infty \left(1-1\wedge\left|1-\frac{3}{n²s²}\right|^n\right)\textnormal{d}s
\end{gather*}
In the penultimate step, one has to notice that, for small s, Doob's inequality gives negative lower bounds for $1-\mathbb{P}\left(\sup_{0\le t<\frac{1}{n}} |\xi_t|²> s\right)$, which have to be improved by the trivial estimate $0\le 1-\mathbb{P}\left(\sup_{0\le t<\frac{1}{n}} |\xi_t|²> s\right)\le 1$.
Now
\begin{equation*}
 \lim_{n\to\infty} \mathbb{E}\sup_{0\le k\le n-1, 0\le t<\frac{1}{n}} \left|\xi_{\frac{k}{n}T+t}-\xi_{\frac{k}{n}T}\right|² =
 \lim_{n\to\infty} \int_0^\infty \left(1-1\wedge\left|1-\frac{3}{n²s²}\right|^n\right)\textnormal{d}s = 0
\end{equation*}
can be shown via the dominated convergence theorem: For the pointwise convergence of the integrand to 0, it is sufficient to show $\lim_{x\to\infty}\left(1-\frac{c}{x²}\right)^x=1$ or, equivalently, \newline $\lim_{x\to\infty}\frac{1}{x^{-1}}\ln\left(1-\frac{c}{x²}\right)=0$ for all $c>0$. This follows from de l'Hospital's theorem because $\frac{\frac{\textnormal{d}}{\textnormal{d}x}\ln\left(1-\frac{c}{x²}\right)}{\frac{\textnormal{d}}{\textnormal{d}x}x^{-1}}=\frac{2c}{\frac{c}{x}-x}\rightarrow 0$. Moreover, for $s\ge 1$ the integrand is monotonically decreasing in n, so that $m(s):= \chi_{[0,1]}(s)+\chi_{(1,\infty)}(s)\left(1-1\wedge\left|1-\frac{3}{s²}\right|\right)= \chi_{[0,1]}(s)+\frac{3}{s²}\chi_{(1,\infty)}(s)$ is an integrable majorant.

Note that the lattice points $\frac{k}{n}T$ play no role and
\begin{gather*}
 \mathbb{E}\sup_{|t-s|<\frac{T}{n}} \left\|\xi_t-\xi_s\right\|² \le
 3\mathbb{E}\sup_{|t-s|<\frac{T}{n}} \left(\|\xi_t-\xi_{\lfloor\frac{nt}{T}\rfloor\frac{T}{n}}\|² + \|\xi_s-\xi_{\lfloor\frac{ns}{T}\rfloor\frac{T}{n}}\|² + \|\xi_{\lfloor\frac{nt}{T}\rfloor\frac{T}{n}}-\xi_{\lfloor\frac{ns}{T}\rfloor\frac{T}{n}}\|² \right) \le \\
 9\mathbb{E}\sup_{0\le k\le n-1, 0\le t\le\frac{1}{n}} \left\|\xi_t-\xi_s\right\|² \to 0 \textnormal{ für } n\to\infty
\end{gather*}
holds as well.

The $m²$ components of the general $L²$-martingale $M_t$ can be represented in the form $M^{i,j}_t=\xi^{i,j}_{\langle M^{i,j}\rangle_t}$ as time-transformed Brownian motions $\xi^{i,j}_t$, $\langle M^{i,j}\rangle_t = \int_0^t |(BC_{0,s})_{i,j}|²\textnormal{d}s$ denoting the quadratic variation. For arbitrary $M>0$, consider the decomposition
\begin{equation*}
 \mathbb{E}\sup_{|t-s|<\frac{T}{n}} \left\|M_t-M_s\right\|²=
 \mathbb{E}\left(\dots\chi_{\sup_{0\le r\le T}||BC_{0,r}||²\le M}\right)+
 \mathbb{E}\left(\dots\chi_{\sup_{0\le r\le T}||BC_{0,r}||²> M}\right) .
\end{equation*}
On the set considered in the first summand, $|\langle M^{i,j}\rangle_t-\langle M^{i,j}\rangle_s| \le M|t-s|$ holds, so the time transform shows that the increments of $M^{i,j}_t$ can be controlled by the ones of a Brownian motion and one gets
\begin{gather*}
 \mathbb{E}\left(\dots\chi_{\sup_{0\le r\le T}||BC_{0,r}||²\le M}\right) \le
 \sum_{i,j} \mathbb{E}\left(\sup_{|t-s|<\frac{T}{n}}|\xi^{i,j}_{\langle M^{i,j}\rangle_t}-\xi^{i,j}_{\langle M^{i,j}\rangle_s}|² \chi_{\sup_{0\le r\le T}||BC_{0,r}||²\le M}\right) \le \\
 \sum_{i,j} \mathbb{E}\left(\sup_{|t-s|<\frac{MT}{n}}|\xi^{i,j}_{t}-\xi^{i,j}_{s}|²\right)
 \xrightarrow{n\to\infty} 0 ;
\end{gather*}
for the second summand there is no such control, so we estimate it at the outset by $\mathbb{E}\left(\sup_{|t-s|\le T} \left\|M_t-M_s\right\|²\chi_{\sup_{0\le r\le T}||BC_{0,r}||²> M}\right)$. According to Doob's inequality, one knows at least
\begin{equation*}
 \mathbb{E}\left(\sup_{|t-s|\le T} \left\|M_t-M_s\right\|²\right) \le
 4\mathbb{E}\left(\sup_{0\le t\le T} \left\|M_t\right\|²\right) \le
 4\mathbb{E}\left(\left\|M_T\right\|²\right) <\infty
\end{equation*}
and since the sets $\{\sup_{0\le r\le T}||BC_{0,r}||²> M\}$ descend to the empty set for $M\to\infty$ because of the continuous paths of $C_{0,t}$,
\begin{equation*}
\lim_{M\to\infty}\mathbb{E} \left(\sup_{|t-s|\le T} \left\|M_t-M_s\right\|²\chi_{\sup_{0\le r\le T}||BC_{0,r}||²> M}\right) = 0
\end{equation*}
follows.

The proof of part a) is therewith complete; in case b), a similar telescopic sum ansatz as in \eqref{Teleskop} yields
\begin{equation}\label{Teleskop2}
\begin{split}
 & g_{n,T}(t)=
 1+\sum_{k=0}^{\lfloor\frac{nt}{T}\rfloor-1}\left(g_{n,T}\left(\frac{k+1}{n}T\right)-g_{n,T}\left(\frac{k}{n}T\right)\right)+
 g_{n,T}(t)-g_{n,T}\left(\left\lfloor\frac{nt}{T}\right\rfloor\frac{T}{n}\right) \\
 =& 1+\sum_{k=0}^{\lfloor\frac{nt}{T}\rfloor-1}\left(F\left(\frac{k}{n}T,\frac{k+1}{n}T\right)-1\right) g_{n,T}\left(\frac{k}{n}T\right)+\left(F\left(\left\lfloor\frac{nt}{T}\right\rfloor\frac{T}{n},t\right)-1\right)
 g_{n,T}\left(\left\lfloor\frac{nt}{T}\right\rfloor\frac{T}{n}\right) \\
 =& 1+\sum_{k=0}^{\lfloor\frac{nt}{T}\rfloor-1}\left(\int_{\frac{k}{n}T}^{\frac{k+1}{n}T}Ag_{n,T}(s)\textnormal{d}s
 +\int_{\frac{k}{n}T}^{\frac{k+1}{n}T}Bg_{n,T}(s)\textnormal{d}\xi_s
 +o_k\left(\frac{1}{n}\right)g_{n,T}\left(\frac{k}{n}T\right)\right)+ \\
 &+\int_{\left\lfloor\frac{nt}{T}\right\rfloor\frac{T}{n}}^t Ag_{n,T}(s)\textnormal{d}s
 +\int_{\left\lfloor\frac{nt}{T}\right\rfloor\frac{T}{n}}^t Bg_{n,T}(s)\textnormal{d}\xi_s
 +o\left(\frac{1}{n}\right)g_{n,T}\left(\left\lfloor\frac{nt}{T}\right\rfloor\frac{T}{n}\right) \\
 =& 1+\int_0^t Ag_{n,T}(s)\textnormal{d}s + \int_0^t Bg_{n,T}(s)\textnormal{d}\xi_s
 +\sum_{k=0}^{\lfloor\frac{nt}{T}\rfloor-1}o_k\left(\frac{1}{n}\right)g_{n,T}\left(\frac{k}{n}T\right)+
 o\left(\frac{1}{n}\right)g_{n,T}\left(\left\lfloor\frac{nt}{T}\right\rfloor\frac{T}{n}\right)
\end{split}
\end{equation}
and one can continue along the lines of the preceding part.
\end{proof}

\begin{cor}\label{stotrotter1}
Formula \eqref{Produkt} holds true, the limit being uniform in the quadratic mean.
\end{cor}

\begin{proof}
 $F(s,t):=e^{(t-s)A}B_{s,t}$ satisfies part b) of the theorem: (ii) folllows from \eqref{det'} and \eqref{stoch'} from the Itô formula one knows
\begin{align*}
  e^{(t-s)A}B_{s,t} &=
  1+\int_s^te^{(\tau-s)A}AB_{s,\tau}\mathrm{d}\tau+\int_s^te^{(\tau-s)A}B_{s,\tau}B\textnormal{d} \xi_\tau \\
  &=1+\int_s^tAe^{(\tau-s)A}B_{s,\tau}\mathrm{d}\tau+\int_s^te^{(\tau-s)A}BB_{s,\tau}\textnormal{d} \xi_\tau \\
  &=1+\int_s^tAe^{(\tau-s)A}B_{s,\tau}\mathrm{d}\tau+\int_s^tBe^{(\tau-s)A}B_{s,\tau}\textnormal{d} \xi_\tau+
  \int_s^t(e^{(\tau-s)A}B-Be^{(\tau-s)A})B_{s,\tau}\textnormal{d} \xi_\tau
\end{align*}
with
\begin{gather*}
 \mathbb{E}\|\frac{1}{t-s}\int_s^t(e^{(\tau-s)A}B-Be^{(\tau-s)A})B_{s,\tau}\textnormal{d}\xi_\tau\|² = \\
 \frac{1}{(t-s)²}\int_s^t\mathbb{E}\|(e^{(\tau-s)A}B-Be^{(\tau-s)A})B_{s,\tau}\|²\textnormal{d}\tau
 \to 0 \ \textnormal{für} \ t\to s
\end{gather*}
because
\begin{gather*}
 \frac{1}{t-s}\mathbb{E}\|(e^{(t-s)A}B-Be^{(t-s)A})B_{s,t}\|² \le
 (t-s)\left\|\frac{e^{(t-s)A}-1}{t-s}B-B\frac{e^{(t-s)A}-1}{t-s}\right\|²\mathbb{E}\|B_{s,t}\|² \\
 \xrightarrow{t\to s}0\cdot\|AB-BA\|²\cdot1=0 .
\end{gather*}
This implies the convergence of $g_{n,T}$. Since $B_{0,t}$ is an $L^2$-martingale, the part of the proof starting with \eqref{f73} tells us that even $\mathbb{E}\sup_{|t-s|\le\frac Tn}\|B_{s,t}\|²\xrightarrow{n\to\infty}0$, so, according to remark \ref{rem:Vorfaktor}, the convergence remains true without the prefactors of $g_{n,T}$ in the form \eqref{Produkt}.
\end{proof}

Theorem \ref{Satz1} also enables us to split off only a part of the deterministic term:

\begin{cor}\label{stotrotter2}
 Let $A=A_1+A_2$. Then $C_{0,t}=\lim_{n\rightarrow\infty}(\prod_{k=0}^{n-1}e^{\frac{t}{n}A_1}D_{\frac{k}{n}t,\frac{k+1}{n}t})$ where $D_{s,t}$ denotes the solution flow of $\textnormal{d}X_t=A_2X_t\textnormal{d}t+BX_t\textnormal{d}\xi_t$.
\end{cor}

\begin{proof}
 Since
 \begin{equation*}
  \begin{split}
   & e^{(t-s)A_1}D_{s,t} =
   1 + \int_s^te^{(\tau-s)A_1}A_1 D_{s,\tau}\mathrm{d}\tau + \int_s^te^{(\tau-s)A_2}A_2 D_{s,\tau}\mathrm{d}\tau
   + \int_s^te^{(\tau-s)A}BD_{s,\tau}\textnormal{d} \xi_\tau \\
   =& 1+\int_s^tAe^{(\tau-s)A}D_{s,\tau}\mathrm{d}\tau+\int_s^tBe^{(\tau-s)A}D_{s,\tau}\textnormal{d} \xi_\tau+
  \int_s^t(e^{(\tau-s)A}B-Be^{(\tau-s)A})D_{s,\tau}\textnormal{d} \xi_\tau ,
  \end{split}
 \end{equation*}
 the proof works as in corollary \ref{stotrotter1}.
\end{proof}

Part a) of the theorem obviously comprises the Euler-Mayurama formula, the notion of convergence being slightly improved w.r.t. \cite{Kloeden}:

\begin{cor}
 $\prod_{k=0}^{n-1} \left(1+ \frac{t}{n}A + (\xi_{\frac{k+1}{n}t}-\xi_{\frac{k}{n}t})B\right) \xrightarrow{n\to\infty}
 C_{0,t}$.
\end{cor}

Moreover, one can approximate the factors in corollary \ref{stotrotter1} by their first-order approximation:

\begin{cor}
 $\prod_{k=0}^{n-1} \left(1+\frac{t}{n}A\right)\left(1+(\xi_{\frac{k+1}{n}t}-\xi_{\frac{k}{n}t})B\right)\xrightarrow{n\to\infty} C_{0,t}$.
\end{cor}

\begin{proof}
 $F(s,t):=(1+(t-s)A)(1+(\xi_t-\xi_s)B)=1+(t-s)A+(\xi_t-\xi_s)B+(t-s)(\xi_t-\xi_s)AB$ satisfies part a) of the theorem.
\end{proof}

As already the one-dimensional example $\textnormal{d}X_t=(1+1)X_t\textnormal{d}\xi_t$ shows, the stochastic part can in  general not be further decomposed: One has $X_t=X_0e^{2\xi_t-2t}$, but the solution of $\textnormal{d}Y_t=Y_t\textnormal{d}\xi_t$ w.r.t. the same initial condition is $Y_t=X_0e^{\xi_t-\frac{t}{2}}$, so the solution flow is $B_{s,t}=e^{(\xi_t-\xi_s)-\frac{(t-s)}{2}}$ and $\prod_{k=0}^{n-1}\left(B_{\frac{k}{n},\frac{k+1}{n}}B_{\frac{k}{n},\frac{k+1}{n}}\right)X_0=e^{2\xi_1-1}X_0\ne X_1$.

\section{Equations with bounded coefficients}
The presented proofs remain true almost without change if one reads \eqref{Gleichung} as an equation for a process with values in a Banach space E, A and B as bounded operators on E, $\|\cdot\|$ as operator norm and the integrals as strong Bochner integrals or, respectively, their stochastic analogue. Only two parts have to be revised:

1. The operator norm is no longer equivalent with a norm generated by a scalar product, as it was used in \eqref{beschr}, so the application of the triangle inequality can no longer be circumvented. But there is another way to improve the estimate \eqref{falsch}: Also by using the submulitplicativity of the norm, one destroys important information about the dynamics of the process - one takes into account the worst case that the oscillation in a time interval $[\frac{k}{n}T, \frac{k+1}{n}T]$ amplifies the preceding one, although they typically cancel. One can start from the representation \eqref{Teleskop} and do already this step via the Gronwall inequality, which maybe does a better job in taking into account the described behaviour: Along the lines of the last section, one gets
\begin{align*}
 \sup_{n\in\mathbb{N}}\mathbb{E}\|f_{n,T}(t)\|² \le&
 4+4\sup_{n\in\mathbb{N}}\mathbb{E}\|\int_0^{\lfloor\frac{nt}{T}\rfloor\frac{T}{n}}Af_{n,T}(s)\textnormal{d}s\|²
  +4\sup_{n\in\mathbb{N}}\mathbb{E}\|\int_0^{\lfloor\frac{nt}{T}\rfloor\frac{T}{n}}Bf_{n,T}(s)\textnormal{d}\xi_s\|² \\
  &+4\sup_{n\in\mathbb{N}}\mathbb{E}\|\sum_{k=0}^{\lfloor\frac{nt}{T}\rfloor-1}o_k\left(\frac{1}{n}\right)f_{n,T}\left(\frac{k}{n}T\right)\|² \\
 \le & 4+4(T\|A\|²+\|B\|²)\int_0^t \sup_{n\in\mathbb{N}}\mathbb{E}\|f_{n,T}(s)\|² \textnormal{d}s
 +o(1)\sup_{n\in\mathbb{N}}\mathbb{E}\|f_{n,T}(t)\|² ;
\end{align*}
in the last summand, a weakening of \eqref{Rest} was used. If N is chosen so large that $o(1) \le \frac{1}{2}$ for $n\ge N$, then
\begin{equation*}
 \frac{1}{2}\sup_{n\ge N}\mathbb{E}\|f_{n,T}(t)\|² \le 4+4(T\|A\|²+\|B\|²)\int_0^t \sup_{n\ge N}\mathbb{E}\|f_{n,T}(s)\|² \textnormal{d}s
\end{equation*}
follows, so $\sup_{t\in [0,T], n\ge N}\mathbb{E}\|f_{n,T}(t)\|² \le 8 e^{8T(T\|A\|²+\|B\|²)}$, which is sufficient since we are interested in $n\to\infty$.

2. The sum appearing in the proof of $\lim_{n\to\infty}\mathbb{E}\sup_{0\le t\le T} \left\| \int_{\lfloor\frac{nt}{T}\rfloor\frac{T}{n}}^t BC_{0,s}\textnormal{d}\xi_s\right\|^2=0$ would now be infinite and not obviously convergent. The claim can still be reduced to its one-dimensional form, namely via the integrator $\xi_t$: For adapted, E-valued processes of the form $a_t(\omega)=\sum_{i=0}^{k-1} a_{\frac{i}{k}T}(\omega)\chi_{\left(\frac{i}{k}T, \frac{i+1}{k}T\right]}(t)$ and bounded by $C>0$, for $n>k$ (that is to say, if no more than one of the points $\frac{i}{k}T$ lies between t and $t+\frac{T}{n}$) either
\begin{equation*}
 \left\| \int_s^t a_r\textnormal{d}\xi_r\right\|^2 =
 \left\| a_{\lfloor\frac{ks}{T}\rfloor\frac{T}{k}} (\xi_t-\xi_s) \right\|^2 \le C² |\xi_t-\xi_s|²
\end{equation*}
(if no one lies in between) or
\begin{equation*}
\begin{split}
 \left\| \int_s^t a_r\textnormal{d}\xi_r\right\|^2 &\le
 2\left\| \int_s^{\lceil\frac{ks}{T}\rceil\frac{T}{k}} a_r\textnormal{d}\xi_r\right\|^2 +
 2\left\| \int_{\lceil\frac{ks}{T}\rceil\frac{T}{k}}^t a_r\textnormal{d}\xi_r\right\|^2 \le
 2C²[(\xi_{\lceil\frac{ks}{T}\rceil\frac{T}{k}}-\xi_s)² + (\xi_t-\xi_{\lceil\frac{ks}{T}\rceil\frac{T}{k}})²] \\
 &\le 4C²\sup_{|t-s|\le \frac{T}{n}}(\xi_t-\xi_s)²
\end{split}
\end{equation*}
(if one lies in between), in any case
\begin{equation*}
 \begin{split}
  \lim_{n\to\infty} \mathbb{E}\sup_{|t-s|\le \frac{T}{n}} \left\| \int_s^t a_r\textnormal{d}\xi_r\right\|^2 \le
  \lim_{n\to\infty} 4C²\mathbb{E}\sup_{|t-s|\le \frac{T}{n}}(\xi_t-\xi_s)² = 0 
 \end{split}
\end{equation*}
holds.
As such processes $a_t$ are dense in $L²(\Omega,E)$, it only remains to show the continuity of the limit in $a_t$. This can be done via Doob's inequality::
\begin{equation*}
\begin{split}
 &\lim_{n\to\infty} \mathbb{E}\sup_{|t-s|\le \frac{T}{n}} \left\| \int_s^t (a_r-\tilde a_r)\textnormal{d}\xi_r\right\|^2  \\
 \le & \lim_{n\to\infty} \mathbb{E}\sup_{|t-s|\le \frac{T}{n}} 2\left(\left\| \int_0^s (a_r-\tilde a_r)\textnormal{d}\xi_r\right\|^2 +
 \left\| \int_0^t (a_r-\tilde a_r)\textnormal{d}\xi_r\right\|^2 \right) \\
 \le &4\lim_{n\to\infty} \mathbb{E}\sup_{0\le t\le T} \left\| \int_0^t (a_r-\tilde a_r)\textnormal{d}\xi_r\right\|^2
 \le 16 \mathbb{E} \left\| \int_0^T (a_r-\tilde a_r)\textnormal{d}\xi_r\right\|^2 \\
 \le &16 \mathbb{E} \int_0^T \|a_r-\tilde a_r\|²\textnormal{d}r
\end{split}
\end{equation*}

\section{Equations with unbounded coefficients}

As outlined in chapter \ref{ch:Loesung}, the stochastic part of an equation like \eqref{linsdg} with unbounded operators should not be considered separately, but rather in connection with a deterministic correction term in the form \eqref{stochTeil}. Therefore a product formula analogous to \eqref{Produkt} can be most easily found for the equation \eqref{stoschroedinger}, which explicitely contains this corrector term.

The proof of the result formulated below shows in particular the weak solvability of \eqref{stoschroedinger2} via the product ansatz \eqref{prostosch}, which is at the same time a more explicit representation of the solution than in \cite{Holevo}. We use ideas from \cite{Holevo}, but our proof is independent. Also we have to formulate ad hoc conditions, but they appear to be rather natural: Roughly speaking, we require the existence of an invariant subspace which the partial solution flows $H_t$ and $A_{s,t}$ do not drive out of the domain of each other's generator. Such conditions appear in deterministic semigroup theory as well, the only difference being that a well-established terminology permits more handsome formulations.

\begin{thm}\label{thm:Produkt}
 Consider the stochastic Schrödinger equation
 \begin{equation}\label{stoschroedinger2}
  \begin{split}
   \textnormal{d}\psi_t = (-iH-\frac{1}{2}A²)\psi_t\textnormal{d}t + A\psi_t\textnormal{d}\xi_t
  \end{split}
 \end{equation}
 with self-adjoint operators H and A and the partial solution flows $H_t:=e^{itH}$ and $A_{s,t}:=e^{(\xi_t-\xi_s)A-(t-s)A²}$. Suppose that there exists a dense subset $\mathcal M\subset\mathfrak{H}$ such that for $\psi\in\mathcal M$ and sufficiently large n all products appearing in \eqref{prostosch} (also the empty one) are in $Dom(A²)$, $A_{t,t+\frac{T}{n}}\prod\ldots\in Dom(H)$ for all t and $H_{-\frac{T}{n}}\prod\ldots\in Dom(A²)$, and another dense set $\mathcal{N}\subset\mathfrak{H}$ such tha for all $\phi\in\mathcal{N}$ both $\lim_{t\to0}AH_t\phi=A\phi$ and $\lim_{t\to0}A²H_t\phi=A²\phi$ holds. Then the equation is weakly solvable and, for all $T>0$ and $\psi\in\mathfrak H$, its solution flow $C_{0,t}$ is given by
 \begin{equation}\label{prostosch}
  C_{0,t}\psi=\lim_{n\to\infty} \prod_{k=0}^{\left\lfloor\frac{nt}{T}\right\rfloor-1} H_{\frac{T}{n}} A_{\frac{k}{n}T,\frac{k+1}{n}T}\psi
 \end{equation}
with weak limit in $L^2(\Omega,\mathfrak{H})$.

If, in addition, the conditions from theorem \ref{starkeLoes} about regular solvability are satisfied and $\mathcal M\cap D(C)$ (C denoting the reference operator used therein) is dense in $\mathfrak{H}$, then the convergence is even strong.
\end{thm}

\begin{proof}
 Choose $\psi\in\mathcal M$ and let $\tilde\psi$ be a product appearing in \eqref{prostosch}. The conditions make sure that, for $s\ge 0, t\in\left(s,s+\frac{1}{n}\right)$, the Itô formula can be applied to $H_{t-s}A_{s,t}\psi$ and, doing this, one gets
 \begin{equation*}
 \begin{split}
  H_{t-s}A_{s,t}\tilde\psi =& \psi-i\int_s^t H_{\tau-s}HA_{s,\tau}\tilde\psi\mathrm d\tau + \int_s^t H_{\tau-s}A_{s,\tau}A\tilde\psi\mathrm d\xi_\tau
  -\frac{1}{2}\int_s^t H_{\tau-s}A_{s,\tau}A²\tilde\psi\mathrm d\tau \\
  =& \tilde\psi+\int_s^t\left(-iH-\frac{1}{2}H_{\tau-s}A²H_{s-\tau}\right)H_{\tau-s}A_{s,\tau}\tilde\psi\mathrm d\tau
  +\int_s^t(H_{\tau-s}AH_{s-\tau})H_{\tau-s}A_{s,\tau}\tilde\psi\mathrm d\xi_\tau .
 \end{split}
 \end{equation*}
 With the ansatz pursued in \eqref{Teleskop2} for $g_{n,T}(t):= H_{t-\left\lfloor\frac{nt}{T}\right\rfloor\frac{T}{n}}A_{\left\lfloor\frac{nt}{T}\right\rfloor\frac{T}{n},t} \prod_{k=0}^{\left\lfloor\frac{nt}{T}\right\rfloor-1} H_{\frac{T}{n}} B_{\frac{k}{n}T,\frac{k+1}{n}T}$ one gets
 \begin{equation}\label{2oj}
 \begin{split}
  g_{n,T}(t)\psi=& \psi+\sum_{k=0}^{\left\lfloor\frac{nt}{T}\right\rfloor-1}
 \left(H_{\frac{T}{n}}A_{\frac{k}{n}T,\frac{k+1}{n}T}-id\right) g_{n,T}\left(\frac{k}{n}T\right)\psi \\
  &+\left(H_{t-\left\lfloor\frac{nt}{T}\right\rfloor\frac{T}{n}}A_{\left\lfloor\frac{nt}{T}\right\rfloor\frac{T}{n},t}-id\right) g_{n,T}\left(\left\lfloor\frac{nt}{T}\right\rfloor\frac{T}{n},t\right)\psi \\
  =&-i\int_0^t H g_{n,T}(s)\mathrm ds
  +\int_0^t H_{s-\left\lfloor\frac{ns}{T}\right\rfloor\frac{T}{n}}AH_{\left\lfloor\frac{ns}{T}\right\rfloor\frac{T}{n}-s}g_{n,T}(s) \mathrm d\xi_s \\
  &-\frac{1}{2} \int_0^t H_{s-\left\lfloor\frac{ns}{T}\right\rfloor\frac{T}{n}}A² H_{\left\lfloor\frac{ns}{T}\right\rfloor\frac{T}{n}-s} g_{n,T}(s) \mathrm ds .
 \end{split}
 \end{equation}
Ths is a linear equation for $g_{n,T}$ with self-adjoint coefficients
\begin{equation*}
\tilde A^{(n,t)}:=H_{t-\left\lfloor\frac{nt}{T}\right\rfloor\frac{T}{n}}AH_{\left\lfloor\frac{nt}{T}\right\rfloor\frac{T}{n}-t} \textnormal{ und }
\tilde{A^2}^{(n,t)}:=H_{t-\left\lfloor\frac{nt}{T}\right\rfloor\frac{T}{n}}A^2H_{\left\lfloor\frac{nt}{T}\right\rfloor\frac{T}{n}-t} ,
\end{equation*}
which, since $0\le s-\left\lfloor\frac{ns}{T}\right\rfloor\frac{T}{n}\le\frac{T}{n}$, converge uniformly in t to their counterparts in \eqref{stoschroedinger2} on $\mathcal{N}$ for $n\to\infty$. Thus, we are in a similar situation as in the proof of theorem \ref{thm:Holevo} (after having arrived at \eqref{linsdgapp}) and, like there, can pass over to the weak limit:

Since $H_t$ and $A_{s,t}$ preserve the $L²(\Omega,\mathfrak H)$-norm,
\begin{equation}\label{iso}
\mathbb{E}\|g_{n,T}(t)\psi\|²=\|\psi\|²
\end{equation}
holds, so the terms $g_{n,T}(\cdot)\psi$, viewed as mappings from a time interval $[0,T]$ to $L^2(\Omega,\mathfrak H)$, are uniformly bounded for any $\psi\in\mathfrak H$.
We would like to show equicontinuity and apply the Arzela-Ascoli theorem. Since we have no clue concerning the limiting behaviour of expressions like
$\mathbb{E}\|\tilde A^{(n,t)}g_{n,T}(t)\psi\|²$, we pass over to weak convergence in $L²(\Omega,\mathfrak{H})\cong L²(\Omega,\mathbb{C})\otimes\mathfrak{H}$, i.e. we choose $\xi\in L²(\Omega,\mathbb{C})$ and $\phi\in\mathcal{N}$ and consider the complex-valued functions $f^n(t):=\mathbb{E}\overline\xi\langle\phi,g_{n,T}(t)\psi\rangle$. Due to the boundedness of $g_{n,T}(t)\psi$, we only have to consider $\xi$ and $\phi$ from dense subsets of the unit spheres, therefore we restrict to $\phi\in\mathcal{N}$; since $L²(\Omega,\mathbb C)$ and $\mathfrak H$ are separable (the first statement follows from the Wiener chaos decomposition), even an appropriate sequence $(\xi_l\phi_l)_{l\in\mathbb N}$ suffices. The $f^n(t)$ are uniformly bounded according to the Cauchy-Schwarz inequality and \eqref{iso} and equicontinuous according to
\begin{equation*}
\begin{split}
 |f^n(t)-f^n(s)|² \le& \mathbb{E}|\xi|²\cdot\mathbb{E}|\langle\phi,g_{n,T}(t)\psi-g_{n,T}(s)\psi\rangle|² \\
 =&\mathbb{E}|\int_s^t\langle\phi,\tilde A^{(n,\tau)}g_{n,T}(\tau)\psi\rangle\mathrm{d}\xi_\tau - \int_s^t\langle\phi, (iH+\frac{1}{2}\tilde{A^2}^{(n,\tau)})g_{n,T}\psi\rangle\mathrm{d}\tau|² \\
 \le& 2\mathbb{E}|\int_s^t\langle \tilde A^{(n,\tau)*}\phi,g_{n,T}(\tau)\psi\rangle\mathrm{d}\xi_\tau|²
 +2\mathbb{E}|\int_s^t\langle (iH+\tilde{A^2}^{(n,\tau)*})\phi, g_{n,T}(\tau)\psi\rangle\mathrm{d}\tau|² \\
 \le& 2\int_s^t\mathbb{E}|\langle \tilde A^{(n,\tau)}\phi,g_{n,T}(\tau)\psi\rangle|²\mathrm{d}\tau
 +2|t-s|\int_s^t\mathbb{E}|\langle (iH-\tilde{A^2}^{(n,\tau)})\phi, g_{n,T}(\tau)\psi\rangle|²\mathrm{d}\tau \\
 \le& 2|t-s|\sup_{\tau\in[s,t]}\|\tilde A^{(n,\tau)}\phi\|²\sup_{\tau\in[s,t]}\mathbb{E}\|g_{n,T}(\tau)\psi\|² \\
 &+2|t-s|²\sup_{\tau\in[s,t]}\|(iH-\tilde{A^2}^{(n,\tau)})\phi\|²\sup_{\tau\in[s,t]}\mathbb{E}\|g_{n,T}(\tau)\psi\|² \\
 \le& 2|t-s|\sup_{\tau\in[s,t]}\|\tilde A^{(n,\tau)}\phi\|²
 +2|t-s|²\sup_{\tau\in[s,t]}\|(iH-\tilde{A^2}^{(n,\tau)})\phi\|² \\
 \le& C_{\phi,T}|t-s|\|\psi\|² .
\end{split}
\end{equation*}
Therefore, there exists a subsequence $(n_k^1)_{k\in\mathbb N}$ of $(n)_{n\in\mathbb N}$ such that $\mathbb{E}\overline\xi_1\langle\phi_1,g_{n_k^1,T}(\cdot)\psi\rangle$ converges uniformly in t, a subsequence $(n_k^2)$ of $(n_k^1)$ such that also $\mathbb{E}\overline\xi_2\langle\phi_2,g_{n_k^2,T}(\cdot)\psi\rangle$ converges and so on. For the diagonal sequence $(n_k^k)=:(n_k)$ obviously $\mathbb{E}\overline\xi_j\langle\phi_j,g_{n_k^k,T}(\cdot)\psi\rangle$ converges for all $j\in\mathbb N$ uniformly in t; since
\begin{equation*}
 L²(\Omega,\mathfrak{H})\ni\eta\mapsto\lim_{k\to\infty}\mathbb{E}\langle\eta,g_{n_k^k,T}(t)\psi\rangle\in\mathbb{R}
\end{equation*}
is a continuous functional, there exists a process $\psi_t$ with
\begin{equation*}
 \lim_{k\to\infty}\mathbb{E}\overline\xi\langle\phi,g_{n_k^k,T}(t)\psi\rangle = \mathbb{E}\overline\xi\langle\phi,\psi_t\rangle
\end{equation*}
for all $\xi,\phi$ uniformly in t.

In order to show that this is the wanted solution, the convergence of the remaining summands in the ``weak form'' of \eqref{2oj}, i.e. in
\begin{align*}
 \mathbb{E}\overline\xi\langle\phi,g_{n_k,T}(t)\psi\rangle =& \mathbb{E}\overline\xi\langle\phi,\psi\rangle
 + \mathbb{E}\overline\xi\langle\phi,\int_0^t \tilde A^{(n_k,s)}g_{n_k,T}(s)\psi\mathrm{d}\xi_s\rangle \\
 &- \mathbb{E}\overline\xi\langle\phi,\int_0^t(iH-\frac{1}{2}\tilde{A^2}^{(n_k,s)}g_{n_k,T}(s)\psi\mathrm{d}s\rangle \textnormal{ resp. }
\end{align*}
\begin{equation}\label{945}
\begin{split}
 \mathbb{E}\overline\xi\langle\phi,g_{n_k,T}(t)\psi\rangle =& \mathbb{E}\overline\xi\langle\phi,\psi\rangle
 + \mathbb{E}\overline\xi\int_0^t\langle \tilde A^{(n_k,s)}\phi,g_{n_k,T}(s)\psi\rangle\mathrm{d}\xi_s \\
 &+\mathbb{E}\overline\xi\int_0^t\langle (iH-\frac 12\tilde{A^2}^{(n_k,s)})\phi,g_{n_k,T}(s)\psi\rangle\mathrm{d}s
\end{split}
\end{equation}
(strong Bochner integrals and their stochastic analogue allow for dragging bounded operators, in particular scalar products, inside),
to their counterparts in the weak version of \eqref{stoschroedinger2} has to be checked.
To begin with,
\begin{align*}
 &|\mathbb{E}\overline\xi\int_0^t\langle (iH-\frac 12\tilde{A^2}^{(n_k,s)})\phi,g_{n_k,T}(s)\psi\rangle\mathrm{d}s
 -\mathbb{E}\overline\xi\int_0^t\langle (iH-\frac 12A^2)\phi,\psi_s\rangle\mathrm{d}s| \\
 =&|\int_0^t\mathbb{E}\overline\xi\langle (iH-\frac 12\tilde{A^2}^{(n_k,s)})\phi,g_{n_k,T}(s)\psi\rangle\mathrm{d}s
 -\int_0^t\mathbb{E}\overline\xi\langle (iH-\frac 12A^2)\phi,\psi_s\rangle\mathrm{d}s| \\
 \le&|\int_0^t\mathbb{E}\overline\xi\langle \frac 12(\tilde{A^2}^{(n_k,s)}\phi-A^2\phi),g_{n_k,T}(s)\psi\rangle\mathrm{d}s|
 +|\int_0^t\mathbb{E}\overline\xi\langle (iH-\frac 12A^2)\phi,(g_{n_k,T}(s)\psi-\psi_s)\rangle\mathrm{d}s| \\
 \le& (\mathbb{E}|\xi|^2)^{\frac{1}{2}} \sup_{s\in[0,t]}\|\frac 12(\tilde{A^2}^{(n_k,s)}\phi-A^2\phi)\| \int_0^t\left(\mathbb{E}\|g_{n_k,T}(s)\psi\|^2\right)^{\frac{1}{2}} \mathrm ds +\hdots
 \to0;
\end{align*}
for the summand in the middle of \eqref{945} we restrict $\xi$ WLOG to the iterated stochastic integrals $I_l(t)$ specified in lemma \ref{lemm:Chaos}. This restriction to a total subset is justified by the fact that, according to
\begin{align*}
 \mathbb{E}|\int_0^t\langle \tilde A^{(n_k,s)}\phi,g_{n_k,T}(s)\psi\rangle\mathrm{d}\xi_s|^2
 =&\int_0^t\mathbb{E}|\langle \tilde A^{(n_k,s)}\phi,g_{n_k,T}(s)\psi\rangle|²\mathrm{d}s \\
 \le& \sup_{s\in[0,T]}\|\tilde A^{(n_k,s)}\phi\|^2 \int_0^t\mathbb{E}\|g_{n_k,T}(s)\psi\|²\mathrm{d}s
 \le CT\|\psi\|^2 ,
\end{align*}
the sequences of stochastic integrals which are checked for convergence are bounded.
Then, similar as above and in addition using the Itô formula, the independence of the Wiener processes and the fact that, according to their recursive definition, all $I_l(t)$ are $L²$-martingales and have mean value 0, we conclude
\begin{align*}
 &|\mathbb{E}\overline{I_{l+1}(t)}\int_0^t\langle \tilde A^{(n_k,s)}\phi,g_{n_k,T}(s)\psi\rangle\mathrm{d}\xi_s
 - \mathbb{E}\overline{I_{l+1}(t)}\int_0^t\langle A\phi,\psi_s\rangle\mathrm{d}\xi_s| \\
 =& |\mathbb{E}\overline{\int_0^t a(s)I_l(s)\mathrm{d}\xi_s}\int_0^t \langle \tilde A^{(n_k,s)}\phi,g_{n_k,T}(s)\psi\rangle\mathrm{d}\xi_s
 - \mathbb{E}\overline{\int_0^t a(s)I_l(s)\mathrm{d}\xi_s}\int_0^t \langle A\phi,\psi_s\rangle\mathrm{d}\xi_s| \\
 =& |\mathbb{E}\int_0^t \overline{a(s)I_l(s)} \langle \tilde A^{(n_k,s)}\phi,g_{n_k,T}(s)\psi\rangle\mathrm{d}s
 - \mathbb{E}\int_0^t \overline{a(s)I_l(s)} \langle A\phi,\psi_s\rangle\mathrm{d}s| \\
 =& |\int_0^t\mathbb{E} \overline{a(s)I_l(s)} \langle (\tilde A^{(n_k,s)}\phi-A\phi,g_{n_k,T}(s)\psi\rangle\mathrm{d}s|
 + |\int_0^t\mathbb{E} \overline{a(s)I_l(s)} \langle A\phi,(g_{n_k,T}(s)\psi-\psi_s)\rangle\mathrm{d}s| \\
 \le& (\mathbb{E}|I_l(s)|^2)^\frac{1}{2}\sup_{s\in[0,T]}\|\tilde A^{(n_k,s)}\phi - A\phi\| \int_0^t\mathbb{E}\|g_{n_k,T}(s)\psi\|^2\mathrm{d}s + \ldots \to0 .
\end{align*}

Therefore, the constructed process $\psi_t$ solves equation \eqref{stoschroedinger2} (and agrees with the solution from \ref{thm:Holevo} insofar as the solution is unique). The same method also allows to start with any subsequence $g_{\tilde n_k,T}(t)\psi$ and extract a sub-subsequence converging in the same manner weakly to $\psi_t$, so even $g_{n,T}(t)\psi \rightharpoonup \psi_t$ holds true. 

Compared to $f_{n,T}(t)\psi:=\prod_{k=0}^{\left\lfloor\frac{nt}{T}\right\rfloor-1} H_{\frac{T}{n}} A_{\frac{k}{n}T,\frac{k+1}{n}T}\psi$, the product representation \eqref{prostosch}, $g_{n,T}$ contains some additional prefactors, but they do not matter because
\begin{align*}
 &|\mathbb{E}\overline\xi\langle\phi, f_{n,T}(t)\psi-g_{n,T}(t)\psi \rangle|^2
 = |\mathbb{E}\overline\xi\langle\phi,
 (H_{t-\left\lfloor\frac{nt}{T}\right\rfloor\frac{T}{n}}A_{\left\lfloor\frac{nt}{T}\right\rfloor\frac{T}{n},t}-1)
 \prod_{k=0}^{\left\lfloor\frac{nt}{T}\right\rfloor-1} H_{\frac{T}{n}} B_{\frac{k}{n}T,\frac{k+1}{n}T}\psi \rangle|^2 \\
 =& |\mathbb{E}\overline\xi\langle (A_{\left\lfloor\frac{nt}{T}\right\rfloor\frac{T}{n},t}^*
 H_{t-\left\lfloor\frac{nt}{T}\right\rfloor\frac{T}{n}}^*-1)\phi,
 \prod_{k=0}^{\left\lfloor\frac{nt}{T}\right\rfloor-1} H_{\frac{T}{n}} B_{\frac{k}{n}T,\frac{k+1}{n}T}\psi \rangle|^2 \\
 \le& \mathbb{E}|\xi|^2 \mathbb{E}\|(A_{\left\lfloor\frac{nt}{T}\right\rfloor\frac{T}{n},t}
 H_{-t+\left\lfloor\frac{nt}{T}\right\rfloor\frac{T}{n}}-1)\phi\|^2
 \mathbb{E}\|\prod_{k=0}^{\left\lfloor\frac{nt}{T}\right\rfloor-1} H_{\frac{T}{n}} B_{\frac{k}{n}T,\frac{k+1}{n}T}\psi\|^2 \\
 =& \mathbb{E}|\xi|^2 \mathbb{E}\|(A_{\left\lfloor\frac{nt}{T}\right\rfloor\frac{T}{n},t}
 H_{-t+\left\lfloor\frac{nt}{T}\right\rfloor\frac{T}{n}}-1)\phi\|^2
 \mathbb{E}\|\psi\|^2 \to0 ,
\end{align*}
since $\left\lfloor\frac{nt}{T}\right\rfloor\frac{T}{n}\to t$ and $A_{s,t}$ and $H_t$, regarded separately, are strongly continuous and uniformly bounded because of $\mathbb{E}\|A_{s,t}\psi\|^2=\mathbb{E}\|H_t\psi\|^2=\|\psi\|^2$ (for $A_{s,t}$ see section \ref{ch:Loesung}.\ref{sec:Einop}).

Since $\lim_{n\to\infty}\mathbb{E}\|f_{n,T}(t)\psi\|²=\|\psi\|²$ holds true and, under the conditions of theorem \ref{thm:Martingal},  $\psi\in\mathcal{M}\cap D(C)$ also satisfies $\mathbb{E}\|\psi_t\|^2=\|\psi\|^2$, in this case the convergence $f_{n,T}(t)\to C_{0,t}$ is even strong. Since $C_{0,t}$ as well as all operators $f_{n,T}$ have $L^2(\Omega,\mathfrak{H})$-norm 1, the strong convergence on the dense subset $\mathcal M\cap D(C)$ implies the convergence on the whole of $\mathfrak{H}$.
\end{proof}

Certainly it would be desirable to prove uniform convergence in the quadrativ mean as in theorem \ref{Satz1}, however, the way to get there is less obvious in this case than in the ansatz using the Gronwall inequality.

In order to get an idea how one can find $\mathcal M$ and $\mathcal N$ in concrete cases, let us again have a look at Diosi's equation
\begin{equation}\label{373}
 \textnormal d\psi_t = -\frac{i}{2}\triangle\psi_t + x\psi_t\textnormal d\xi_t
 - \frac{1}{2}x^2\psi_t\textnormal dt .
\end{equation}
We try the maximally possible choice
\begin{equation*}
\mathcal M=\mathcal N=D(\triangle)\cap D(x²)=\{\phi\in H^2(\mathbb{R})\mid \cdot^2\phi(\cdot)\in L^2\}=:\mathcal D .
\end{equation*}
If $\phi\in\mathcal D$, then obviously $A_{s,t}\phi\in\mathcal D$, since
\begin{equation*}
 A_{s,t}\phi(x)=e^{x(\xi_t-\xi_s)-x^2(t-s)}\phi(x)
\end{equation*}
and multicplication with a Gaussian preserves decay and smoothness properties. $H_t\phi\in D(\triangle)$ ist well-known and instead of $H_t\phi\in D(x^2)$ it is easier to show $\widehat{H_t\phi}\in H^2(\mathbb{R})$ (Sobolev space): This follows from the fact that $\widehat{H_t\phi}$ solves the Fourier transformed Schrödinger equation $\dot{\hat{\psi_t}}(\xi) = -\frac{i}{2}\xi^2\hat{\psi_t}(\xi)$ and is therefore given by $\widehat{H_t\phi} = e^{-\frac{i}{2}\cdot^2t}\hat\phi$. It remains to be shown that $\mathcal D$ is also a suitable choice for $\mathcal N$; according to the Plancherel formula, one has to check $\widehat{xH_t\phi(x)}\xrightarrow{L^2}\widehat{x\phi(x)}$ and $\widehat{x^2H_t\phi(x)}\xrightarrow{L^2}\widehat{x^2\phi(x)}$ for $t\to0$ or, equivalently, $\frac{d}{d\xi}\widehat{H_t\phi}(\xi)\xrightarrow{L^2}\frac{d}{d\xi}\widehat{\phi}(\xi)$ and $\frac{d^2}{d\xi^2}\widehat{H_t\phi}(\xi)\xrightarrow{L^2}\frac{d^2}{d\xi^2}\widehat{\phi}(\xi)$. The formula for $\widehat{H_t\phi}$ tells us
\begin{align*}
 &\frac{d}{d\xi}\widehat{H_t\phi}(\xi) = e^{-\frac{i}{2}\xi^2t}(-it\xi\hat\phi(\xi)+\frac{d}{d\xi}\hat\phi(\xi))
 \textnormal{ and} \\
 &\frac{d^2}{d\xi^2}\widehat{H_t\phi}(\xi) = e^{-\frac{i}{2}\xi^2t}
 [-(t^2\xi^2+it)\hat\phi(\xi)-2it\xi\frac{d}{d\xi}\hat\phi(\xi)+\frac{d^2}{d\xi^2}\hat\phi(\xi)]
\end{align*}
and the assertions follow from the dominated convergence theorem since $\phi\in H^2$ and $\cdot\phi(\cdot)\in H^1$ imply that  $\cdot^2\hat\phi$ and $\cdot\frac{d}{d\xi}\hat\phi(\cdot)$ are in $L^2$.

The fact that $\mathcal D$ agrees with the domain \eqref{Domref} of the mentionned reference operator implies strong convergence in the product formula.

The application in chapter \ref{ch:Beispiel} requires the following straightforward generalization of theorem \ref{thm:Produkt}:
\begin{cor}\label{cor:Produkt}
 Consider the stochastic Schrödinger equation
 \begin{equation}\label{stoschroedinger3}
  \begin{split}
   \textnormal{d}\psi_t = -iH\psi_t\textnormal{d}t + \sum_{j=1}^m(A_j\psi_t\textnormal{d}\xi_t^j-\frac{1}{2}A_j^2\psi_t\textnormal dt)
  \end{split}
 \end{equation}
 with independent Wiener processes $\xi_t^j$, self-adjoint operator H and commuting self-adjoint operators $A_j$ and the partial solution flows $H_t:=e^{itH}$ and
\begin{equation*}
A_{s,t}:=e^{\sum_j[(\xi_t^j-\xi_s^j)A_j-(t-s)A_j^2]} = \prod_j e^{(\xi_t^j-\xi_s^j)A_j-(t-s)A_j^2} .
\end{equation*}
Suppose that there exists a dense subset $\mathcal M\subset\mathfrak{H}$ such that for $\psi\in\mathcal M$ and sufficiently large n all products appearing in \eqref{prostosch2} (also the empty one) are in $\bigcap_j Dom(A_j^2)$, $A_{t,t+\frac{T}{n}}\prod\ldots\in Dom(H)$ for all t and $H_{-\frac{T}{n}}\prod\ldots\in\bigcap Dom(A_j^2)$, and another dense set $\mathcal{N}\subset\mathfrak{H}$such that for all $\phi\in\mathcal{N}$ and j both $\lim_{t\to0}A_jH_t\phi=A_j\phi$ and $\lim_{t\to0}A_j^2H_t\phi=A_j^2\phi$ holds. Then the equation is weakly solvable and, for all $T>0$ and $\psi\in\mathfrak H$, its solution flow $C_{0,t}$ is given by
 \begin{equation}\label{prostosch2}
  C_{0,t}\psi=\lim_{n\to\infty} \prod_{k=0}^{\left\lfloor\frac{nt}{T}\right\rfloor-1} H_{\frac{T}{n}} A_{\frac{k}{n}T,\frac{k+1}{n}T}\psi
 \end{equation}
with weak limit in $L^2(\Omega,\mathfrak{H})$.

If, in addition, the conditions of theorem \ref{starkeLoes} about regular solvability are satisfied and $\mathcal M\cap D(C)$ (C denoting the reference operator used therein) is dense in $\mathfrak{H}$, then the convergence is even strong.
\end{cor}

\begin{proof}
 As in the proof of theorem \ref{thm:Produkt} one gets by means of the Itô formula
\begin{align*}
 &H_{t-s}A_{s,t}\tilde\psi = \\
 =& \psi-i\int_s^t H_{\tau-s}HA_{s,\tau}\tilde\psi\mathrm d\tau + \sum_j(\int_s^t H_{\tau-s}A_{s,\tau}A_j\tilde\psi\mathrm d\xi_\tau^j
  -\frac{1}{2}\int_s^t H_{\tau-s}A_{s,\tau}A_j^2\tilde\psi\mathrm d\tau) \\
  =& \tilde\psi+\int_s^t-iHH_{\tau-s}A_{s,\tau}\tilde\psi\mathrm d\tau \\
  &+\sum_j(\int_s^t(H_{\tau-s}AH_{s-\tau})H_{\tau-s}A_{s,\tau}\tilde\psi\mathrm d\xi_\tau
  -\frac 12 \int_s^t (H_{\tau-s}A²H_{s-\tau})H_{\tau-s}A_{s,\tau}\tilde\psi\mathrm d\tau .
\end{align*}
(Here it is important that the $A_j$ commute - otherwise, terms like
\begin{equation*}
\prod_{j=1}^{k-1} e^{(\xi_t^j-\xi_s^j)A_j-(t-s)A_j^2} A_k \prod_{j=k}^m e^{(\xi_t^j-\xi_s^j)A_j-(t-s)A_j^2}
\end{equation*}
would appear and would not agree with $A_kA_{s,t}$ because $A_k$ could not be put in front. In order to deal with such cases, one would not only have to require compatibility conditions between $H$ and $A_j$, but likewise for the $A_j$ among themselves.)

Via the same telescopic sum ansatz one gets the equation
\begin{align*}
 g_{n,T}(t)\psi=&-i\int_0^t H g_{n,T}(s)\mathrm ds
  +\sum_j\Big(\int_0^t H_{s-\left\lfloor\frac{ns}{T}\right\rfloor\frac{T}{n}}A_jH_{\left\lfloor\frac{ns}{T}\right\rfloor\frac{T}{n}-s}g_{n,T}(s) \mathrm d\xi_s \\
  &-\frac{1}{2} \int_0^t H_{s-\left\lfloor\frac{ns}{T}\right\rfloor\frac{T}{n}}A_j^2 H_{\left\lfloor\frac{ns}{T}\right\rfloor\frac{T}{n}-s} g_{n,T}(s) \mathrm ds\Big) ,
\end{align*}
the limit of which can be performed as in the previous proof. The $\xi$ needed for testing weak convergence are still provided by lemma \ref{lemm:Chaos}.
\end{proof}
Equation \eqref{QMUPL-Kollaps}, the multi-dimensional version of Diosi's equation \eqref{373}, is covered by this corollary in the same way as \eqref{373} by \ref{thm:Produkt}.

\chapter{Example: Quantum mechanical collapse models - GRW and QMUPL}\label{ch:Beispiel}

\section{The models and the question}

In the end, I would like to present the two quantum mechanical collapse models which made me think about the product formula. Both models complement the conventional Schrödinger evolution $\phi_t=e^{-\frac{i}{h}tH}\phi$ by a - mathematically well-defined - mechanism leading to a random spatial localization of the wavefunction und thus solving the ``measurement problem'' which is often illustrated by ``Schrödinger's cat''. $\phi_t$ is replaced by a stochastic process (we will not change the notation, though) on a probability space $(\Omega,\mathfrak{F},\mathbb{P})$. Arbitrary Schrödinger operators H are admitted, though, for simplicity, we will restrict to the case of a single free particle and omit purely physical constants, so $\psi\in L²(\mathbb{R}³,\mathbb{C})$ and $H=-\frac{1}{2}\triangle$.

\subsection{GRW}
In the Ghirardi-Rimini-Weber model (shortly GRW, see \cite{Ghirardi1}, \cite{Bassi} and, for a mathematically rigorous treatment, \cite{Tumulka}), the Schrödinger evolution is interrupted at the random jump times $\tau_1, \tau_2 \dots$ of a Poisson process with intensity $\mu$ (so, in paricular, $\mathbb{E}\tau_k=\frac{k}{\mu}$) by a collapse mechanism - namely, the multiplication by a Gaussian ``hitting function'' with spread $\frac{1}{\sqrt\alpha}$. For mathematical simplicity and - see later - in accordance with our needs - we replace the $\tau_k$ by deterministic times $\frac{k}{\mu}$ Then the mechanism can be described as follows:

We start with an initial wavefunction $\phi_0\in L²(\mathbb{R}^3,\mathbb{C})$, construct a process $(\phi_{\frac{k}{\mu}})$ of wavefunctions on the space $\Omega:=(\mathbb{R}^3)^{\{\frac{1}{\mu},\frac{2}{\mu},\dots\}}$, the coordinate projections of which will be denoted $Y_1, Y_2, \dots$, recursively via
\begin{equation}\label{GRW}
 \begin{split}
 & \psi_n(Y_1, \dots, Y_n, x) := \left(\frac{\alpha}{\pi}\right)^{\frac{1}{4}} e^{-\frac{\alpha}{2}|x-Y_n|²} e^{-\frac{i}{\mu}H} \phi_{n-1}(Y_1, \dots, Y_{n-1}, x) \\
 & \phi_n:=\frac{\psi_n}{\|\psi_n\|_{L²(\mathbb{R}^3)}}
 \end{split}
\end{equation}
and define
\begin{equation}\label{flash}
\begin{split}
 \mathbb{P}_{\alpha,\mu}(Y_n\in A\mid Y_1, \dots, Y_{n-1})
 :=& \int_A \|\psi_n(Y_1,\dots,Y_{n-1},y,\cdot)\|_{L^2(\mathbb{R}^3)}^2 \textnormal dy \\
 =& \int_A \sqrt{\frac{\alpha}{\pi}}\int e^{-\alpha|x-y|²}|(e^{-\frac{i}{\mu}H}\phi_{n-1})(Y, x)|²\textnormal{d}x\textnormal{d}y ,
\end{split}
\end{equation}
which reads explicitely as
\begin{equation}\label{flash1}
\begin{split}
 \mathbb{P}_{\alpha,\mu}(Y_1\in A_1) =& \int_{A_1}\|\psi_1(y_1,\cdot)\|_{L^2(\mathbb{R}^3)}^2 \textnormal dy_1 \\
 =& \int_{A_1} \sqrt{\frac{\alpha}{\pi}}\int e^{-\alpha|x-y_1|²}|(e^{-\frac{i}{\mu}H}\phi_0)(x)|²\textnormal{d}x\textnormal{d}y_1 ,
\end{split}
 \end{equation}
\begin{equation*}
\begin{split}
 &\mathbb{P}_{\alpha,\mu}(Y_1\in A_1, Y_2 \in A_2) \\
 =& \int_{\{Y_1\in A_1\}} \mathbb{P}_{\alpha,\mu}(Y_2\in A_2\mid Y_1)\textnormal{d}\mathbb{P}_{\alpha,\mu} \\
 =& \int_{\{Y_1\in A_1\}} \int_{A_2}\sqrt{\frac{\alpha}{\pi}}\int e^{-\alpha|x-y_2|^2}
 |(e^{-\frac{i}{\mu}H}\frac{\psi_1}{\|\psi_1\|})(Y_1,x)|^2 \textnormal dx\textnormal dy_2 \textnormal d\mathbb{P}_{GRW}^{Y_1} \\
 =& \int_{A_1} \int_{A_2}\sqrt{\frac{\alpha}{\pi}}\int e^{-\alpha|x-y_2|^2} |(e^{-\frac{i}{\mu}H}\psi_1)(y_1,x)|^2
 \textnormal dx\textnormal dy_2 \frac{1}{\|\psi_1(y_1,\cdot)\|^2} \|\psi_1(y_1,\cdot)\|^2 \textnormal dy_1 \\
 =& \int_{A_1} \int_{A_2}\frac{\alpha}{\pi}\int e^{-\alpha|x-y_2|^2} |(e^{-\frac{i}{\mu}H}e^{-\frac{\alpha}{2}|\cdot-y_1|^2}e^{-\frac{i}{\mu}H}\phi_0)(x)|^2
 \textnormal dx\textnormal dy_2 \textnormal dy_1
\end{split}
\end{equation*}
and, inductively,
\begin{equation}\label{flashes}
\begin{split}
 & \mathbb{P}_{\alpha,\mu}(Y_1\in A_1, \dots, Y_n \in A_n) \\
 =& \int_{A_1\times\cdots\times A_n}\sqrt{\frac{\alpha}{\pi}}^n \int |(e^{-\frac{\alpha}{2}|\cdot-y_n|^2}e^{-\frac{i}{\mu}H} \cdots
 e^{-\frac{\alpha}{2}|\cdot-y_1|^2}e^{-\frac{i}{\mu}H}\phi_0)(x)|^2 \textnormal dx\textnormal d(y_1, \dots, y_n)
\end{split}
\end{equation}
and yields a consistent family that has a projective limit on the whole of $\Omega$. The GRW process is obtained by canonically identifying $\Omega$ with a subset of $C([0,\infty),\mathbb{R}^3)$, pushing $\mathbb{P}_{\alpha,\mu}$ forward to this space via the corresponding inclusion and extending $(\phi_t)$ according to $\phi_{\frac{k}{\mu}+s}:=e^{-isH}\phi_{\frac{k}{\mu}}$ between the collapse times.

\subsection{QMUPL}
In the QMUPL (Quantum Mechanics with Universal Position Localization) model, tracing back to Diosi (\cite{Diosi1}, \cite{Diosi2}), a time-continuous collapse of an initial wavefunction $\phi_0\in L²$ is constructed by the aid of the stochastic Schrödinger equation
\begin{equation}\label{QMUPL-Kollaps}
 \begin{split}
 &d\psi_t(x)=-iH\psi_t\textnormal dt + \sqrt{\lambda}\psi_t(x)x\cdot \textnormal d\xi_t - \frac{\lambda}{2}|x|²\psi_t(x) \textnormal dt \\
 &\psi_0=\phi_0
 \end{split}
\end{equation}
with a three-dimensional Wiener process $(\xi_t)$ on a filtered probability space $(\Omega, \mathfrak{F, (F_t)}, \mathbb{Q})$ and an intensity parameter $\lambda>0$: The ``physical'' collapse process is $\phi_t:=\frac{\psi_t}{\|\psi_t\|_{L²(\mathbb{R}^3)}}$, weighted by a new measure defined by
\begin{equation}\label{cooking}
\mathbb{P}_\lambda(A):= \mathbb{E}_\mathbb{Q}(\chi_A\|\psi_t\|²) \textnormal{ für } A\in\mathfrak{F}_t .
\end{equation}
The martingale property of $\|\psi_t\|²$ (theorem \ref{thm:Martingal}) ensures that this is indeed a consistent definition of a measure on $\Omega$. One physical justification of this model lies in the fact that the isolated collapse process (i.e. \eqref{QMUPL-Kollaps}, \eqref{cooking} with H=0) is a continuous extension of the GRW collapse process (see the after next paragraph).

\subsection{The question}
In \cite{Diosi1}, Diosi has claimed that the QMUPL model can be obtained as a scaling limit from GRW by increasing the collapse frequency and, at the same time, the spread of the hitting function, thus weakening the effect of a single collapse and keeping the overall effect constant. The appropriate scaling of the GRW process for $\alpha\to0$ was found to be 
\begin{equation}\label{Skalierung}
\mu=\frac{2\lambda}{\alpha}
\end{equation}
with constant $\lambda>0$. Diosis proof lacks mathematical rigor and only considers the level of the Lindblad equations for the statistical operator $\rho_t:=\mathbb{E}(|\psi_t\rangle\langle\psi_t|)$, which are (in ``position representation'')
\begin{equation*}
 \frac{d}{dt}\rho_t(x,y) = -\frac{i}{\hbar}[H,\rho_t](x,y) -\mu\left(1-e^{-\frac{\alpha}{4}(x-y)²}\right)\rho_t(x,y)
 \end{equation*}
for GRW and
\begin{equation*}
 \frac{d}{dt}\rho_t(x,y)= -\frac{i}{\hbar}[H,\rho_t](x,y) -\frac{\lambda}{2}(x-y)²\rho_t(x,y)
\end{equation*}
for QMUPL. For $\alpha\to0$, the GRW Lindblad term can be linearized in $\alpha$ and in fact one gets
\begin{equation*}
 \mu\left(1-e^{-\frac{\alpha}{4}(x-y)²}\right)\approx \mu\left(1-(1-\frac{\alpha}{4}(x-y)²\right) = \frac{\lambda}{2}(x-y)² .
\end{equation*}
We have paved the way for proving the following result directly on the level of the wavefunction:
\begin{thm}\label{thm:Limes}
 For $\alpha\to0$, the finite dimensional distributions $\mathbb{P}_{\alpha,\frac{2\lambda}{\alpha}}\circ(\phi_{t_1},\dots,\phi_{t_n})^{-1}$ of the scaled GRW processes converge weakly to $\mathbb{P}_{\lambda}\circ(\phi_{t_1},\dots,\phi_{t_n})^{-1}$, which are those of the QMUPL process.
\end{thm}
It seems plausible that, apart form technical complications, random GRW collapse times would not change anything. Moreover, it is worth mentionning that the only purpose of the random collapse times is to avoid an arbitrary assignment and, in the case of several particles, an arbitrary coupling of the collapse times for the single particles. This arbirariness is removed by the continuum limit anyway.

\section{Equivalence of the collapse processes}
The comparison of GRW and QMUPL is particularly elucidating if one neglects the Schrödinger evolution (i.e. puts H=0). Then, according to \eqref{GRW} and \eqref{flashes},
\begin{equation}\label{871}
 \phi_{\frac{n}{\mu}}(x)=c e^{-\frac{\alpha}{2}[|x-Y_1|²+\cdots+|x-Y_n|²]}\phi_0(x)
\end{equation}
(c is the Y-dependent normalization factor) with
\begin{equation}\label{872}
 \mathbb{P}_{\alpha,\mu}(Y_1\in A_1,\dots,Y_n\in A_n)
 = \sqrt{\frac{\alpha}{\pi}}^n \int_{A_1\times\cdots\times A_n} \int e^{-\alpha[(x-y_1)²+\dots+(x-y_n)²]}|\phi_0(x)|²
 \textnormal dx\textnormal d(y_1, \dots, y_n) .
\end{equation}
The solution of the collapse part of \eqref{QMUPL-Kollaps} is
\begin{equation*}
 \psi_t(x) = e^{\sqrt{\lambda} x\cdot\xi_t - \lambda |x|²t}\phi_0(x) ,
\end{equation*}
and can be rewritten as
\begin{equation*}
\begin{split}
 \psi_{\frac{n}{\mu}}(x) =& e^{\sum_{k=1}^n[\sqrt{\lambda}x\cdot(\xi_{\frac{k}{\mu}}-\xi_{\frac{k-1}{\mu}}) -
 \frac{\lambda}{\mu}|x|²]} \phi_0(x)
 = e^{-\frac{\lambda}{\mu}\sum_{k=1}^n [|x-\frac{\mu}{2\sqrt{\lambda}}(\xi_{\frac{k}{\mu}}-\xi_{\frac{k-1}{\mu}})|² 
 + \frac{\mu^2}{4\lambda}|\xi_{\frac{k}{\mu}}-\xi_{\frac{k-1}{\mu}}|²]} \phi_0(x)
\end{split}
\end{equation*}
or, in normalized form and taking into account the scaling \eqref{Skalierung},
\begin{equation}\label{873}
 \phi_{\frac{n}{\mu}}(x)
 = c e^{-\frac{\alpha}{2}\sum_{k=1}^n [|x-\frac{\mu}{2\sqrt{\lambda}}(\xi_{\frac{k}{\mu}}-\xi_{\frac{k-1}{\mu}})|²]} \phi_0(x) .
\end{equation}
The quantities $Z_k:=\frac{\mu}{2\sqrt{\lambda}}(\xi_{\frac{k}{\mu}}-\xi_{\frac{k-1}{\mu}})$ formally correspond to the $Y_k$ from GRW and, by the aid of \eqref{cooking}, their distribution under $\mathbb P_\lambda$ can be determined if one takes into account that they have independent centered normal distributions with $\mathbb{V}(Z_k)=\frac{\mu}{4\lambda}=\frac{1}{2\alpha}$ under $\mathbb{Q}$:
\begin{equation}\label{874}
 \begin{split}
  &\mathbb{P}_\lambda(Z_1\in A_1, \dots, Z_n\in A_n)
  = \mathbb{E}_\mathbb{Q}(\chi_{\{Z_1\in A_1, \dots, Z_n\in A_n\}}\|\psi_\frac{n}{\mu}\|^2) \\
  =& \mathbb{E}_\mathbb{Q}(\chi_{\{Z_1\in A_1, \dots, Z_n\in A_n\}} \int
  e^{\sum_{k=1}^n[2\sqrt{\lambda}x\cdot(\xi_{\frac{k}{\mu}}-\xi_{\frac{k-1}{\mu}}) -
  \frac{2\lambda}{\mu}|x|²]} |\phi_0|^2(x) \textnormal dx ) \\
  =& \mathbb{E}_\mathbb{Q}(\chi_{\{Z_1\in A_1, \dots, Z_n\in A_n\}} \int
  e^{\sum_{k=1}^n[\frac{4\lambda}{\mu}x\cdot Z_k -
  \frac{2\lambda}{\mu}|x|²]} |\phi_0|^2(x) \textnormal dx ) \\
  =& \int_{A_1\times\dots\times A_n} \int e^{\sum_{k=1}^n[2\alpha x\cdot z_k - \alpha|x|²]} |\phi_0|^2(x) \textnormal dx
  \sqrt{\frac{\alpha}{\pi}}^n e^{-\alpha\sum_{k=1}^n |z_k|^2} \textnormal d(z_1,\dots,z_n) \\
  =& \sqrt{\frac{\alpha}{\pi}}^n \int_{A_1\times\dots\times A_n} \int e^{-\alpha\sum_{k=1}^n|x-z_k|^2} |\phi_0|^2(x) \textnormal dx
  \textnormal d(z_1,\dots,z_n)
 \end{split}
\end{equation}
\eqref{873} and \eqref{874} agree with \eqref{871} and \eqref{872}. Thus, there is no qualitative difference between GRW and QMUPL processes with parameters $(\alpha,\frac{2\lambda}{\alpha})$ and $\lambda$  - the first can simply be obtained by restricting the latter to the appropriate discrete instants of time.

\section{QMUPL as continuum limit of GRW}
The circumstance that the GRW collapse is a discretization of the QMUPL one suggests to read the iterative definition of the complete GRW model (including the Hamiltonian) as a product formula approximation of the solution of \eqref{QMUPL-Kollaps} - the only question to be settled is how the change of measure \eqref{cooking} fits into this framework. To this end, we not only approximate the solution of the QMUPL equation up to a time T, according to theorem \ref{thm:Produkt}, by
\begin{equation}\label{480}
 \psi_t^\mu:= \prod_{k=0}^{\lfloor\mu t\rfloor-1}e^{\sqrt\lambda (\cdot)\cdot(\xi_\frac{k+1}{\mu}-\xi_\frac{k}{\mu})-\frac{\lambda}{\mu}|\cdot|^2}
 e^{-\frac{i}{\mu}H}\phi_0 ,
\end{equation}
but also the change of measure by
\begin{equation*}
 \mathbb{P}_\lambda^\mu(A):=\mathbb{E}_\mathbb{Q}(\chi_A\|\psi_t^\mu\|^2) \textnormal{ für } A\in\mathfrak{F}_T .
\end{equation*}
Then, with $Z_k:=\frac{\mu}{2\sqrt{\lambda}}(\xi_{\frac{k}{\mu}}-\xi_{\frac{k-1}{\mu}})$,
\begin{equation*}
 \frac{\psi_t^\mu}{\|\psi_t^\mu\|_{L²}}=c\prod_{k=0}^{\lfloor\mu t\rfloor-1} e^{-\frac{\alpha}{2}|\cdot-Z_{k+1}|^2}e^{-\frac{i}{\mu}H}\phi
\end{equation*}
with $\xi$-dependent normalization factor c and, as a generalization of \eqref{874}, one finds that
\begin{equation*}
 \begin{split}
  &\mathbb{P}_\lambda^\mu(Z_1\in A_1,\dots, Z_n\in A_n)
  = \mathbb{E}_\mathbb{Q}(\chi_{\{Z_1\in A_1, \dots, Z_n\in A_n\}}\|\psi_\frac{n}{\mu}^\mu\|^2) \\
  =& \mathbb{E}_\mathbb{Q}(\chi_{\{Z_1\in A_1, \dots, Z_n\in A_n\}} \int |\prod_{k=0}^{n-1}
  e^{\frac{2\lambda}{\mu}(\cdot)\cdot Z_{k+1} - \frac{\lambda}{\mu}|\cdot|²}e^{-\frac{i}{\mu}H}\phi|^2(x) \textnormal dx ) \\
  =& \int_{A_1\times\dots\times A_n} \int |\prod_{k=0}^{n-1}
  e^{\alpha (\cdot)\cdot z_{k+1} - \frac{\alpha}{2}|\cdot|²}e^{-\frac{i}{\mu}H}\phi|^2(x) \textnormal dx
  \sqrt{\frac{\alpha}{\pi}}^n e^{-\alpha\sum_{k=0}^{n-1} |z_{k+1}|^2} \textnormal d(z_1,\dots,z_n) \\
  =& \sqrt{\frac{\alpha}{\pi}}^n\int_{A_1\times\dots\times A_n} \int |\prod_{k=0}^{n-1}
  e^{\alpha (\cdot)\cdot z_{k+1} - \frac{\alpha}{2}|\cdot|² - \frac{\alpha}{2}|z_{k+1}|²}e^{-\frac{i}{\mu}H}\phi|^2(x) \textnormal dx
  \textnormal d(z_1,\dots,z_n) \\
  =& \sqrt{\frac{\alpha}{\pi}}^n\int_{A_1\times\dots\times A_n} \int |\prod_{k=0}^{n-1}
  e^{-\frac{\alpha}{2}|\cdot-z_{k+1}|^2}e^{-\frac{i}{\mu}H}\phi|^2(x) \textnormal dx \textnormal d(z_1,\dots,z_n)
 \end{split}
\end{equation*}
if $\frac{n}{\mu}\le T$ ($(\cdot)$ denotes the argument of the function). This expression agrees with the GRW distribution \eqref{flashes}, i.e. on $\mathfrak{F}_T$ it holds that $\mathbb{P}_{\lambda}^\mu=\mathbb{P}_{\alpha,\mu}$ with $\alpha=\frac{2\lambda}{\mu}$. The question of convergence in theorem \ref{thm:Limes} for $\alpha\to0$ resp. $\mu\to\infty$ can now be settled without further reference to the GRW model by showing
\begin{equation*}
 \lim_{\mu\to\infty}\mathbb{E}(\chi_{\{\phi_{t_1}^\mu\in A_1, \dots, \phi_{t_n}^\mu\in A_n\}}\|\psi_{t_n}^\mu\|_{L^2})
 =\mathbb{E}(\chi_{\{\phi_{t_1}\in A_1, \dots, \phi_{t_n}\in A_n\}}\|\psi_{t_n}\|_{L^2}) ,
\end{equation*}
$\psi_t$ denoting the solution of the QMUPL equation \eqref{QMUPL-Kollaps}. However, this is obvious because, since $\psi_t^\mu\to\psi_t$ in $L^2(\Omega, L^2(\mathbb{R}^3))$ according to corollary \ref{cor:Produkt}, the sequence of the integrands is convergent in probability and uniformly integrable.

Note that the reverse order of the factors in the definition \eqref{480} of $\psi_t^\mu$ w.r.t. theorem \ref{thm:Produkt} plays no role:
\begin{align*}
 &\prod_{k=0}^{\left\lfloor\frac{nt}{T}\right\rfloor-1} H_{\frac{T}{n}} A_{\frac{k}{n}T,\frac{k+1}{n}T}\psi
 -\prod_{k=0}^{\left\lfloor\frac{nt}{T}\right\rfloor-1} A_{\frac{k}{n}T,\frac{k+1}{n}T} H_{\frac Tn}\psi \\
 =& (H_{\frac Tn}-1) A_{\left\lfloor\frac{nt}{T}\right\rfloor-1,\left\lfloor\frac{nt}{T}\right\rfloor} \prod_{k=0}^{\left\lfloor\frac{nt}{T}\right\rfloor-2} H_{\frac{T}{n}} A_{\frac{k}{n}T,\frac{k+1}{n}T}\psi
 +\left(\prod_{k=1}^{\left\lfloor\frac{nt}{T}\right\rfloor-1} A_{\frac{k}{n}T,\frac{k+1}{n}T} H_{\frac Tn}\psi\right)A_{0,\frac Tn}(1-H_{\frac Tn})\psi ;
\end{align*}
for the second summand, the unitarity of $H_t$ and \eqref{isometrie} imply $\mathbb{E}\|\cdot\|^2=\|(1-H_{\frac Tn})\psi\|^2\to0$ for $n\to\infty$; as for the first, we know at least
\begin{align*}
 &|\mathbb{E}\overline\xi\langle \phi, (H_{\frac Tn}-1) A_{\left\lfloor\frac{nt}{T}\right\rfloor-1,\left\lfloor\frac{nt}{T}\right\rfloor} \prod_{k=0}^{\left\lfloor\frac{nt}{T}\right\rfloor-2} H_{\frac{T}{n}} A_{\frac{k}{n}T,\frac{k+1}{n}T}\psi \rangle|^2 \\
 =& |\mathbb{E}\overline\xi\langle \left(\prod_{k=\left\lfloor\frac{nt}{T}\right\rfloor-2}^{0} A_{\frac{k}{n}T,\frac{k+1}{n}T}^* H_{\frac{T}{n}}^*\right)  A_{\left\lfloor\frac{nt}{T}\right\rfloor-1,\left\lfloor\frac{nt}{T}\right\rfloor}^*(H_{\frac Tn}^*-1)\phi, \psi \rangle|^2 \\
 \le& \mathbb{E}|\xi|^2\mathbb{E}\|\left(\prod_{k=\left\lfloor\frac{nt}{T}\right\rfloor-2}^{0} A_{\frac{k}{n}T,\frac{k+1}{n}T} H_{-\frac{T}{n}}\right)
 A_{\left\lfloor\frac{nt}{T}\right\rfloor-1,\left\lfloor\frac{nt}{T}\right\rfloor}(H_{-\frac Tn}-1)\phi\|^2\|\psi\|^2 \\
 =& \mathbb{E}|\xi|^2\|(H_{-\frac Tn}-1)\phi\|^2\|\psi\|^2\to0 .
\end{align*}
Thus, despite the different order, $\psi_t^\mu$ converges weakly and, since $\mathbb{E}\|\psi_t^\mu\|^2=\|\psi\|^2=\mathbb{E}\|\psi_t\|^2$, also strongly to $\psi_t$.

\chapter{Conclusion}

We have found a new access to the solution of equation \eqref{346} by means of a representation which, compared to earlier works, is more explicit and closer to the physical motivation of such equations. Possibilities for technical improvements are obvious and werw mentioned on the spot. The applicability of our product formula to the QMUPL equation subject to the nonlinear modification by a measure change, as described in chapter \ref{ch:Beispiel}, raises the general question to what extent such product formulas remain valid for diffusion processes described by nonlinear equations (in the deterministic theory of nonlinear semigroups, this topic is inquired e.g. in \cite{Kato2}). In particular, since the QMUPL model, as it seems, does not admit adaptions to further physical requirements (such as identical particles, Lorentz invariance, relativistic interaction etc.), an extension of the product formula to models like CSL (\cite{Ghirardi2}), which are considered as more promising in this respect, is suggested.

\backmatter
\bibliographystyle{plain}
\bibliography{Quellen.bib}

\end{document}